\PassOptionsToPackage{linktocpage}{hyperref}
\documentclass[review]{elsarticle}
\usepackage[fleqn]{amsmath}
\usepackage{amssymb,latexsym,amsmath}

\usepackage{overpic}
\usepackage{color}
\usepackage{hyperref}
\usepackage{graphicx}
\usepackage{amsfonts}
\usepackage{subfigure}
\usepackage{epstopdf}
\usepackage{subfigure}
\usepackage{epsfig}
\graphicspath{{../Plots/}}

\makeatletter
\def\ps@pprintTitle{%
   \let\@oddhead\@empty
   \let\@evenhead\@empty
   \let\@oddfoot\@empty
   \let\@evenfoot\@oddfoot
}
\makeatother

\newtheorem{thm}{Theorem}[section]
\newtheorem{defn}{Definition}[section]

\newtheorem{lem}{Lemma}[section]

\newtheorem{exmp}{Example}[section]
\numberwithin{equation}{section}

\begin{document}
\begin{frontmatter}
\title{ Local discontinuous Galerkin method for   distributed-order time and space-fractional
convection-diffusion  and Schr\"{o}dinger type equations}
\author[]{Tarek Aboelenen}
\ead{tarek.aboelenen@aun.edu.eg}




\address{Department of Mathematics, Assiut University, Assiut 71516, Egypt}
\begin{abstract}
 Fractional partial differential equations with distributed-order fractional derivatives describe some important physical phenomena. In this paper,  we propose  a local discontinuous Galerkin (LDG)  method   for   the distributed-order time and Riesz space fractional convection-diffusion  and Schr\"{o}dinger type equations. We prove stability and optimal order of convergence $\mathcal{O}(h^{N+1}+(\Delta t)^{1+\frac{\theta}{2}}+\theta^{2})$ for the distributed-order time and  space-fractional diffusion  and Schr\"{o}dinger type equations, an order of convergence of $\mathcal{O}(h^{N+\frac{1}{2}}+(\Delta t)^{1+\frac{\theta}{2}}+\theta^{2})$ is established for the distributed-order time and Riesz space fractional
convection-diffusion equations  where $\Delta t$,  $h$ and  $\theta$ are  the  step  sizes  in  time,
space and distributed-order variables, respectively.  Finally, the performed numerical experiments confirm the optimal order of convergence.

{\bf Keywords:} \emph{ time distributed order and space-fractional  convection-diffusion equations, time distributed order and space-fractional  Schr\"{o}dinger type equations, local discontinuous Galerkin  method, stability, error estimates.}\\


\end{abstract}

\end{frontmatter}


\section{Introduction}

The distributed-order differential equation can serve  as a natural generalization of
the single-order and the multi-term fractional differential equation \cite{diethelm2009numerical} which arises in many physical and biological applications, for example, the stress behavior of an elastic medium\cite{caputo1967linear}, the
torsional phenomenon of anelastic or dielectric spherical shells and infinite planes \cite{Caputo}, the rheological properties of composite materials \cite{hartley1999fractional,lorenzo2002variable}, dielectric induction and diffusion \cite{caputo2001distributed},  viscoelastic oscillators \cite{atanackovic2005fractional},
distributed order membranes in the ear \cite{naghibolhosseini2015estimation}, and anomalous diffusion \cite{caputo2003diffusion,Sokolov}.  The earliest appearance of the idea on the distributed-order equation may date back to the Caputo's work in 1960s \cite{Caputo},
 which was also stated by Podlubny et al.  \cite{podlubny2013matrix}.  Jiao et al. \cite{JiaoDistributed}  presented a concise and insightful view to understand the usefulness of distributed-order concept in control and signal processing. a more complicated process cannot be described by a single power law and a mixture of power laws leads to a time distributed-order fractional derivative  \cite{meerschaert2011distributed}. Chechkin et al. \cite{chechkin2002retarding} proposed diffusion-like
equations with time and space fractional derivatives of the distributed order for the kinetic description of anomalous diffusion and relaxation phenomena and demonstrated that retarding subdiffusion
and accelerating superdiffusion were governed by distributed-order fractional diffusion equation.
Luchko \cite{luchko2009boundary} investigated some uniqueness and existence results of solutions to boundary value problems
of the generalized distributed-order time-fractional diffusion equation by an appropriate maximum principle. Gorenflo et al. \cite{gorenflo2013fundamental} obtained a representation of the fundamental solution to the Cauchy problem of a distributed-order time-fractional diffusion-wave equation by employing the technique of the Fourier and Laplace transforms and gave the interpretation of the fundamental solution as a probability density function.
Furthermore, they studied waves in a viscoelastic rod of finite
length, where viscoelastic material was described by a constitutive equation of fractional distributed-order type (see Atanackovic et al. \cite{atanackovic2011distributed}).\\
In recent years, developing various numerical algorithms for solving distributed-order and space-fractional equations has received much attention. For the distributed-order time differential equations, Diethelm and Ford \cite{diethelm2009numerical,diethelm2001numerical} presented the numerical methods for solving the distributed-order ordinary differential equations, where the
distributed-order integral was ﬁrstly approximated using the quadrature formula and then
the multi-term fractional differential equations were resulted in, which were ﬁnally reduced
to a system of single-term equations. The idea was followed by Ford and Morgado \cite{ford2012distributed} still
for the distributed-order ordinary differential equations. The matrix approach to the solution
of distributed-order differential equations was introduced by Podlubny et al. \cite{podlubny2013matrix}. For the distributed-order time and space-fractional  equations, Ye et al.
\cite{ye2013numerical} have treated the time distributed-order and space Riesz fractional diffusion on bounded
domains numerically, where the distributed integral was discretized by the mid-point quadrature rule and the time-fractional derivatives in the resultant multi-term fractional diffusion
equation were approximated by the classical L1 formula. Hu et al. \cite{hu2016implicit} investigated an implicit numerical method for the time distributed-order and two-sided space-fractional advection-dispersion equation.  Jacobi collocation method in two successional steps is developed to numerically solve the
multi-dimensional distributed-order generalized Schr\"{o}dinger equations \cite{bhrawy2017numerical}.
To the best of our knowledge, however, the LDG method, which is an important approach to solve partial differential equations  and fractional partial differential equations , has not been considered for the distributed-order time and space-fractional  partial equations.
In this paper, we develop a LDG method to solve the distributed-order time and space-fractional convection-diffusion equations equation
\begin{equation}\label{25n}
\begin{split}
&\mathcal{D}_{t}^{W(\alpha)}{}u+\varepsilon(-\Delta)^{\frac{\beta}{2}}u+\frac{\partial}{\partial x} f(u)=0, \quad x\in \mathbb{R},\,\,t\in(0,T],\\
&u(x,0) = u_{0}(x),\quad x\in \mathbb{R},
\end{split}
\end{equation}
 the nonlinear distributed-order time and space-fractional   Schr\"{o}dinger  equation
\begin{equation}\label{sch1vn}
\begin{split}
&i\mathcal{D}_{t}^{W(\alpha)}{}u- \varepsilon_{1}(-\Delta)^{\frac{\beta}{2}}u+ \varepsilon_{2}f(|u|^{2})u=0,\quad x\in \mathbb{R},\,\,t\in(0,T],\\
&u(x,0) = u_{0}(x),\quad x\in \mathbb{R},
\end{split}
\end{equation}
and the coupled nonlinear distributed-order time and space-fractional    Schr\"{o}dinger equations
\begin{equation}\label{sch2}
\begin{split}
&i\mathcal{D}_{t}^{W(\alpha)}{}u- \varepsilon_{1}(-\Delta)^{\frac{\beta}{2}}u+ \varepsilon_{2}f(|u|^{2},|v|^{2})u=0,\quad x\in \mathbb{R},\,\,t\in(0,T],\\
&i\mathcal{D}_{t}^{W(\alpha)}{}v- \varepsilon_{3}(-\Delta)^{\frac{\beta}{2}}v+\varepsilon_{4} g(|u|^{2},|v|^{2})v=0,\quad x\in \mathbb{R},\,\,t\in(0,T],\\
&u(x,0) = u_{0}(x),\quad x\in \mathbb{R},\\
&v(x,0) =  v_{0}(x),\quad x\in \mathbb{R},
\end{split}
\end{equation}
and homogeneous boundary conditions.  $f(u)$ and $g(u)$ are arbitrary (smooth) nonlinear real functions and  $\varepsilon, \varepsilon_{i}$, $i=1,2,3,4$  are a real constants,  and $\mathcal{D}_{t}^{W(\alpha)}{}u(x,t)$ denotes the distributed order fractional derivative of $u$ in time $t$, given by
\begin{equation}\label{25nzv}
\begin{split}
&\mathcal{D}_{t}^{W(\alpha)}{}u=\int_{0}^{1}W(\alpha){}^{\,\,\,C}_{\,\,\,\,0}\mathcal{D}_{t}^{\alpha}u(x,t)d\alpha,
\end{split}
\end{equation}
where $W(\alpha)$ is the weight function, ${}^{\,\,\,C}_{\,\,\,\,0}\mathcal{D}_{t}^{\alpha}u(x,t)$ $(0 < \alpha < 1)$ is the Caputo fractional derivative of order $\alpha$ with respect to $t$. The  fractional Laplacian $-(-\Delta)^{\frac{\beta}{2}}$, which can be defined using Fourier analysis as  \cite{el2006finite,muslih2010riesz, yang2010numerical}
$$-(-\Delta)^{\frac{\beta}{2}}u(x,t)=\mathcal{F}^{-1}(|\xi|^{\beta}\hat{u}(\xi,t)),$$
where $\mathcal{F}$ is the Fourier transform.\\
The discontinuous Galerkin (DG) method is a class of finite element methods using
discontinuous, piecewise polynomials as the solution and the test spaces in
the spatial direction. There have been various DG methods suggested in the literature to solve diffusion problem, including the method originally proposed by Bassi and Rebay \cite{bassi1997high} for compressible Navier-Stokes
equations, its generalization called the local discontinuous Galerkin (LDG) methods introduced
in \cite{cockburn1989tvb} by Cockburn and Shu and further studied in \cite{cockburn2002approximation, cockburn2005locally}.These DG methods have several attractive properties. It can be easily designed
for any order of accuracy and it has the advantage of greatly facilitates the handling of
complicated geometries and elements of various shapes and types, as well as the treatment
of boundary conditions. And the higher-order of convergence can be achieved without over
many iterations.  For application of the method to fractional problems, Mustapha and McLean \cite{Mustapha2011,mustapha2012uniform,mustapha2013superconvergence,deng2013local}  have developed and analyzed discontinuous Galerkin methods for time fractional diffusion and wave equations.  Xu and Hesthaven \cite{doi:10.1137/130918174}  proposed a LDG method for  fractional convection-diffusion
equations. They proved stability and optimal order of convergence $N+1$
for the fractional diffusion problem when polynomials of degree $N$, and an order of convergence of $N+\frac{1}{2}$ is established for the general fractional convection-diffusion problem with general monotone flux for the nonlinear term. Aboelenen and El-Hawary \cite{cann} proposed a high-order nodal discontinuous Galerkin  method for  a linearized fractional  Cahn-Hilliard equation. They proved stability and optimal order of convergence $N+1$
for the linearized fractional Cahn-Hilliard  problem. A nodal discontinuous Galerkin  method  was   developed to  solve the nonlinear
 Riesz space fractional Schr\"{o}dinger equation and  the strongly coupled nonlinear Riesz space fractional Schr\"{o}dinger equations  \cite{Aboelenen2018428}. They proved, for both problems, $L^{2}$ stability  and optimal order of convergence $O(h^{N+1})$. Aboelenen  \cite{AboelenenDDG} proposed a direct discontinuous Galerkin   (DDG) finite element method for  fractional convection-diffusion and Schr\"{o}dinger type equations. they  proved, for both problems,  $L^{2}$ stability  and  a priori $L^{2}$ error estimates.\\
 This paper is organized as follows. In Section \ref{s1}, we introduce some basic definitions and recall a few central results. We derive the discontinuous Galerkin formulation for the  distributed-order time and Riesz space fractional convection-diffusion equations in Section \ref{s2}. Then  we prove a theoretical result of $L^{2}$ stability  as well as an error estimate  in Section \ref{s3}. In Section \ref{s4}, we  present and analyze a local discontinuous Galerkin method for the nonlinear distributed-order time and Riesz space    fractional Schr\"{o}dinger  type equations. We derive the discontinuous Galerkin formulation for the nonlinear  distributed-order time and Riesz space fractional  Schr\"{o}dinger equation in Section \ref{s41}. Moreover, we prove a theoretical result of $L^{2}$ stability for the nonlinear case in Section \ref{s42} as well as an error estimate for the linear case in Section \ref{s43}. In Section \ref{s44}, we present a local discontinuous Galerkin method for the nonlinear  distributed-order time and Riesz space   fractional coupled nonlinear Schr\"{o}dinger equations and give a theoretical result of $L^{2}$ stability and error estimates. Section \ref{s5} presents some numerical examples to  illustrate the efficiency of the scheme. A few concluding remarks are offered in Section \ref{s6}.

\section{Preliminary definitions}\label{s1}
 We introduce some preliminary definitions of fractional calculus, see, e.g.,\cite{miller1993introduction} and associated functional setting for the subsequent numerical
schemes and theoretical analysis.
 \subsection{Liouville-Caputo fractional calculus}
 The left-sided and right-sided Riemann-Liouville integrals of order $\mu$, when $0 < \mu < 1$, are defined, respectively, as
 \begin{equation}\label{1}
\big({}^{\,\,RL}_{-\infty}\mathcal{I}_{x}^{\mu}f\big)(x)=\frac{1}{\Gamma(\mu)}\int_{-\infty}^{x}
\frac{f(s)ds}{(x-s)^{1-\mu}}, \quad x > -\infty,
\end{equation}
and
\begin{equation}\label{1111}
\big({}^{RL}_{\,\,\,\,x}\mathcal{I}_{\infty}^{\mu}f\big)(x)=\frac{1}{\Gamma(\mu)}\int_{x}^{\infty}\frac{f(s)ds}{(s-x)^{1-\mu}}, \quad x < \infty,
\end{equation}
where $\Gamma$ represents the Euler Gamma function. The corresponding inverse operators, i.e., the left-sided and
right-sided  fractional derivatives of order $\mu$, are then defined based on \eqref{1} and \eqref{1111}, as
\begin{equation}\label{2}
\big({}^{\,\,RL}_{-\infty}\mathcal{D}_{x}^{\mu}f\big)(x)=\frac{d}{dx}\big({}^{\,\,RL}_{-\infty}
\mathcal{I}_{x}^{1-\mu}f\big)(x)=\frac{1}{\Gamma(1-\mu)}
\frac{d}{dx}\int_{-\infty}^{x}\frac{f(s)ds}{(x-s)^{\mu}}, \quad x > -\infty,
\end{equation}
and
\begin{equation}\label{3}
\big({}^{RL}_{\,\,\,\,x}\mathcal{D}_{\infty}^{\mu}f\big)(x)=\frac{-d}{dx}\big(^{RL}_{\,\,\,\,x}\mathcal{I}_{\infty}
^{1-\mu}f\big)(x)=\frac{1}{\Gamma(1-\mu)}
\bigg(\frac{-d}{dx}\bigg)\int_{x}^{\infty}\frac{f(s)ds}{(s-x)^{\mu}}, \quad x < \infty.
\end{equation}
This allows for the definition of the left and right Riemann-Liouville fractional derivatives of order $\mu$ $ (n-1 <\mu<n),\,\,n\in \mathbb{N}$ as
\begin{equation}\label{2}
\big({}^{\,\,RL}_{-\infty}\mathcal{D}_{x}^{\mu}f\big)(x)=\bigg(\frac{d}{dx}\bigg)^{n}\big({}^{\,\,RL}_{-\infty}
\mathcal{I}_{x}^{n-\mu}f\big)(x)=\frac{1}{\Gamma(n-\mu)}
\bigg(\frac{d}{dx}\bigg)^{n}\int_{-\infty}^{x}\frac{f(s)ds}{(x-s)^{-n+1+\mu}}, \quad x > -\infty,
\end{equation}
and
\begin{equation}\label{3}
\big({}^{RL}_{\,\,\,\,x}\mathcal{D}_{\infty}^{\mu}f\big)(x)=\bigg(\frac{-d}{dx}\bigg)^{n}\big(^{RL}_{\,\,\,\,x}
\mathcal{I}_{\infty}
^{n-\mu}f\big)(x)=\frac{1}{\Gamma(n-\mu)}
\bigg(\frac{-d}{dx}\bigg)^{n}\int_{x}^{\infty}\frac{f(s)ds}{(s-x)^{-n+1+\mu}}, \quad x < \infty.
\end{equation}
Furthermore, the corresponding left-sided and right-sided  Caputo derivatives of order $\mu$ $ (n-1 <\mu<n)$ are obtained as
\begin{equation}\label{4}
\big({}^{\,\,\,\,\,\,C}_{-\infty}\mathcal{D}_{x}^{\mu}f\big)(x)=\bigg({}^{RL}_{-\infty}\mathcal{I}_{x}^{n-\mu}
\frac{d^{n}f}{dx^{n}}\bigg)(x)
=\frac{1}
{\Gamma(n-\mu)}\int_{-\infty}^{x}\frac{f^{(n)}(s)ds}{(x-s)^{n-1+\mu}}, \quad x > -\infty,
\end{equation}
and
\begin{equation}\label{5}
\big({}^{C}_{\,x}\mathcal{D}_{\infty}^{\mu}f\big)(x)=(-1)^{n}\bigg({}^{RL}_{\,\,\,\,x}
\mathcal{I}_{\infty}^{n-\mu}\frac{d^{n}f}{dx^{n}}
\bigg)(x)=
\frac{1}{\Gamma(n-\mu)}\int_{x}^{\infty}\frac{(-1)^{n}f^{(n)}(s)ds}{(s-x)}^{n-1+\mu}, \quad x < \infty.
\end{equation}
To carry out the analysis, we introduce the appropriate fractional spaces.
 \begin{defn}(left fractional space \cite{Ervin_variationalformulation}). We define the seminorm
\begin{equation}\label{7}
|u|_{J_{L}^{\mu}(\mathbb{R})}=\big\|{}^{RL}_{\,x_{L}}\mathcal{D}_{x}^{\mu}u\big\|_{L^{2}(\mathbb{R})}.
\end{equation}
and the norm
\begin{equation}\label{7}
\|u\|_{J_{L}^{\mu}(\mathbb{R})}=(|u|_{J_{L}^{\mu}(\mathbb{R})}^{2}+\|u\|_{L^{2}(\mathbb{R})}^{2})^{\frac{1}{2}},
\end{equation}
and let $J_{L}^{\mu}(\mathbb{R})$ denote the closure of $C_{0}^{\infty}(\mathbb{R})$ with respect to $\|.\|_{J_{L}^{\mu}(\mathbb{R})}$.\\
\end{defn}
\begin{defn} (right fractional space \cite{Ervin_variationalformulation}). We define the seminorm
\begin{equation}\label{7}
|u|_{J_{R}^{\mu}(\mathbb{R})}=\big\|{}^{RL}_{\,\,\,\,x}
\mathcal{D}_{x_{R}}^{\mu}u\big\|_{L^{2}(\mathbb{R})},
\end{equation}
and the norm
\begin{equation}\label{7}
\|u\|_{J_{R}^{\mu}(\mathbb{R})}=(|u|_{J_{R}^{\mu}(\mathbb{R})}^{2}+\|u\|_{L^{2}(\mathbb{R})}^{2})^{\frac{1}{2}},
\end{equation}
and let $J_{R}^{\mu}(\mathbb{R})$ denote the closure of $C_{0}^{\infty}(\mathbb{R})$ with respect to $\|.\|_{J_{R}^{\mu}(\mathbb{R})}$.
\end{defn}
\begin{defn} (symmetric fractional space \cite{Ervin_variationalformulation}). We define the seminorm
\begin{equation}\label{7}
\|u\|_{J_{S}^{\mu}(\mathbb{R})}=\big|\big({}^{RL}_{\,x_{L}}\mathcal{D}_{x}^{\mu}u,{}^{RL}_{\,\,\,\,x}
\mathcal{D}_{x_{R}}^{\mu}u\big)_{L^{2}(\mathbb{R})}\big|^{\frac{1}{2}},
\end{equation}
and the norm
\begin{equation}\label{7}
\|u\|_{J_{S}^{\mu}(\mathbb{R})}=\big(|u|_{J_{S}^{\mu}(\mathbb{R})}^{2}+\|u\|_{L^{2}(\mathbb{R})}^{2}\big)^{\frac{1}{2}}.
\end{equation}
and let $J_{S}^{\mu}(\mathbb{R})$ denote the closure of $C_{0}^{\infty}(\mathbb{R})$ with respect to $\|.\|_{J_{S}^{\mu}(\mathbb{R})}$.
\end{defn}
\begin{lem}\label{lk}(see \cite{Ervin_variationalformulation}). For any $0 <s<1$, the fractional integral satisfies the following
property:
\begin{equation}\label{7}
({}^{\,\,RL}_{-\infty}\mathcal{I}_{x}^{s}u,{}^{RL}_{\,\,\,\,x}\mathcal{I}_{\infty}^{s}u)_{\mathbb{R}}=
\cos(s\pi)|u|_{J_{L}^{-s}(\mathbb{R})}^{2}=
\cos(s\pi)|u|_{J_{R}^{-s}(\mathbb{R})}^{2}.
\end{equation}
\end{lem}
Generally, we consider the problem in a bounded domain instead of $\mathbb{R}$. Hence,
we restrict the definition to the domain $\Omega = [a, b]$.
\begin{defn} Define the spaces $J_{R,0}^{\mu}(\Omega),J_{L,0}^{\mu}(\Omega),J_{S,0}^{\mu}(\Omega)$ as the closures of
$C_{0}^{\infty}(\Omega)$ under their respective norms.
\end{defn}
\begin{lem}\label{lg}
(fractional Poincar$\acute{e}$-Friedrichs, \cite{Ervin_variationalformulation}). For $u \in J_{L,0}^{\mu}(\Omega)$ and $\mu \in \mathbb{R}$, we have
\begin{equation}\label{7}
\|u\|_{L^{2}(\Omega)}\leq C |u|_{J_{L,0}^{\mu}(\Omega)},
\end{equation}
and for $u \in J_{R,0}^{\mu}(\Omega)$, we have
\begin{equation}\label{7}
\|u\|_{L^{2}(\Omega)}\leq C |u|_{J_{R,0}^{\mu}(\Omega)}.
\end{equation}
\end{lem}
\begin{lem}\label{lga2} (See \cite{Kilbas:2006:TAF:1137742})
The fractional integration operator $\mathcal{I}^{s}$ is bounded in $L^{2}(\Omega)$:
\begin{equation}\label{7}
\|\mathcal{I}^{s}u\|_{L^{2}(\Omega)}\leq K \|u\|_{L^{2}(\Omega)},
\end{equation}
where $\mathcal{I}^{s}={}^{RL}_{\,x_{L}}\mathcal{I}_{x}^{s}$ (i.e., right-sided Riemann-Liouville integral  of order $s$).
\end{lem}
\begin{lem}\label{lga2} (See \cite{Aboelenen2018428})
The fractional integration operator $\Delta_{-\mu}$ is bounded in $L^{2}(\Omega)$:
\begin{equation}\label{7}
\|\Delta_{-\mu}u\|_{L^{2}(\Omega)}\leq K \|u\|_{L^{2}(\Omega)}.
\end{equation}
\end{lem}
\section{  LDG scheme for  the time distributed-order and space-fractional
  convection-diffusion equation}\label{s2}
Let us consider  the distributed-order time and Riesz space fractional  convection-diffusion  equation.
We first discretize the integral interval $[0,1]$ by the grid $0=\tau_{0}<\tau_{1}<...<\tau_{S}=1$ and take
$\Delta\tau_{j}=\tau_{j}-\tau_{j-1}=\frac{1}{S}=\theta$, $\alpha_{j}=\frac{\tau_{j}+\tau_{j-1}}{2}=\frac{2j-1}{2S}$, $j=1,2,...,S$, $S\in\mathbb{N}$. Then using the mid-point quadrature rule, we obtain
\begin{equation}\label{7zz}
\begin{split}
&\mathcal{D}_{t}^{W(\alpha)}u(x,t)=\sum_{j=1}^{S}W(\alpha_{j}){}^{\,\,\,C}_{\,\,\,\,0}
\mathcal{D}_{t}^{\alpha_{j}}u(x,t)\Delta\tau_{j}
+\mathcal{O}(\theta^{2}),
\end{split}
\end{equation}
where $\theta$ is the step size of the discretization of the numerical integration. Thus the
distributed-order fractional equation \eqref{25n} is now transformed into  multi-term fractional equation.
An approximation to the time fractional derivative \eqref{7zz} can be obtained by simple quadrature formula given as \cite{sun2006fully}.  Let $\Delta t=T/M$ be the time mesh-size, $M$ is a positive integer, $t_{n}=n\Delta t,\,n=0,1,...,M$ be mesh points.
\begin{lem}\label{lga2z} (See \cite{sun2006fully})
Suppose $(0<\alpha<1)$, $y(t)\in\mathcal{C}^{2}[0,t_{n}]$. It holds that
\begin{equation}\label{7}
\begin{split}
&\bigg|\frac{1}{\Gamma(1-\alpha)}\int_{0}^{t_{n}}\frac{y^{'}(s)ds}{(t_{n}-s)^{\alpha}}-
\frac{1}{\lambda}\bigg[a_{0}y(t_{n})-\sum_{l=1}^{n-1}(a_{n-l-1}-a_{n-l})y(t_{l})-a_{n-1}y(0)\bigg]\bigg|\\
&\leq \frac{1}{\Gamma(2-\alpha)}\bigg[\frac{1-\alpha}{12}+\frac{2^{2-\alpha}}{2-\alpha}-(1+2^{-\alpha})\bigg]
\max_{0\leq t\leq t_{n}}|y^{''}(t)|(\Delta t)^{2-\alpha}.
\end{split}
\end{equation}
\end{lem}
For simplicity of the presentation of the proposed method, we introduce the notation
\begin{equation}\label{1bvvzzc}
\begin{split}
&{}^{\,\,\,C}_{\,\,\,\,0}\mathcal{D}_{t_{n}}^{\alpha} y\approx\delta_{t}^{\alpha}y^{n}=\frac{1}{\lambda}\bigg(y^{n}-\sum_{l=1}^{n-1}(a_{n-l-1}-a_{n-l})y^{l}-a_{n-1}y^{0}\bigg).
\end{split}
\end{equation}
From \eqref{25nzv}, \eqref{7zz} and \eqref{1bvvzzc} we obtain
\begin{equation}\label{1bvvzz}
\begin{split}
\mathcal{D}_{t_{n}}^{W(\alpha)}u&\approx\sum_{j=1}^{S}\Delta\tau_{j}
W(\alpha_{j})
{}^{\,\,\,C}_{\,\,\,\,0}\mathcal{D}_{t_{n}}^{\alpha_{j}}u\approx\sum_{j=1}^{S}
\Delta\tau_{j}W(\alpha_{j})
\delta_{t}^{\alpha_{j}}u^{n}\\
&=\sum_{j=1}^{S}
\frac{W(\alpha_{j})\Delta\tau_{j}}{\lambda_{j}}\bigg(u^{n}-\sum_{l=1}^{n-1}(a_{n-l-1}^{\alpha_{j}}-
a_{n-l}^{\alpha_{j}})u^{l}
-a_{n-1}^{\alpha_{j}}u^{0}\bigg),
\end{split}
\end{equation}
where $\lambda_{j}=(\Delta t)^{\alpha_{j}} \Gamma(2-\alpha_{j})$ and $a_{l}^{\alpha_{j}}=(l+1)^{1-\alpha_{j}}-l^{1-\alpha_{j}}$, $0\leq l\leq M-1$.\\
 To obtain a high order discontinuous Galerkin scheme for the space fractional derivative, we rewrite the
fractional derivative as a composite of first order derivatives and a fractional integral
to recover the equation to a low order system. However, for the first order system,
alternating fluxes are used. We introduce  variables $p$, $q$ and $r$ and set
\begin{equation}\label{1a}
\begin{split}
&p=\Delta_{(\beta-2)/2}q,\quad q=\frac{\partial}{\partial x}r,\quad r=\frac{\partial}{\partial x}u,\\
\end{split}
\end{equation}
then, the time distributed order and space-fractional  convection-diffusion problem can be rewritten as
\begin{equation}\label{1b}
\begin{split}
&\mathcal{D}_{t}^{W(\alpha)}u+\frac{\partial}{\partial x}f(u)-p=0,\\
&p=\Delta_{(\beta-2)/2}q,\quad q=\frac{\partial}{\partial x}r,\quad r=\frac{\partial}{\partial x}u.\\
\end{split}
\end{equation}
Now we introduce the broken Sobolev space for any real number $r$
\begin{equation}\label{82aa1}
H^{r}(\Omega)=\{v\in L^{2}(\Omega):\forall k=1,2,....K,v|_{D^{k}}\in H^{r}(D^{k})\}.
\end{equation}
We define the local inner product and $L^{2}(D^{k})$ norm
\begin{equation}\label{82aa1}
(u,v)_{D^{k}}=\int_{D^{k}}uvdx,\quad \|u\|^{2}_{D^{k}}=(u,u)_{D^{k}},
\end{equation}
as well as the global broken inner product and norm
\begin{equation}\label{82aa1}
(u,v)=\sum_{k=1}^{K}(u,v)_{D^{k}},\quad \|u\|^{2}_{L^{2}(\Omega)}=\sum_{k=1}^{K}(u,u)_{D^{k}}.
\end{equation}
We introduce some notation
\begin{equation}\label{82aa1}
u^{\pm}(x_{i})=\lim_{x\rightarrow x_{i}^{\pm}}u(x),\quad\{u\}=\frac{u^{+}+u^{-}}{2},\quad [u]=u^{+}-u^{-}.
\end{equation}
For  simplicity we  discretize  the  computational domain $\Omega$ into $K$ non-overlapping elements,  $D^{k}=[x_{k-\frac{1}{2}},x_{k+\frac{1}{2}}]$, $k = 1,...,K$. Let $u_{h}^{n}, p_{h}^{n}, q_{h}^{n}, r_{h}^{n}\in V_{k}^{N}$ be the approximation of $u(.,t_{n}), p(.,t_{n}),q(.,t_{n}),r(.,t_{n})$ respectively, where  the
approximation space  is defined as
\begin{equation}\label{82aa1}
V_{k}^{N}=\{v:v_{k}\in\mathbb{P}(D^{k}), \, \forall D^{k}\in \Omega\},
\end{equation}
where $\mathbb{P}(D^{k})$ denotes the set of polynomials of degree up to $N$  defined  on  the  element  $D^{k}$.

 We define a fully discrete local discontinuous Galerkin scheme with as follows: find $u_{h}^{n},p_{h}^{n}, q_{h}^{n}, r_{h}^{n}\in V_{k}^{N}$, such that for all test functions $v,\psi,\phi,\eta\in V_{k}^{N}$,
\begin{equation}\label{1dv}
\begin{split}
&\bigg(\sum_{j=1}^{S}W(\alpha_{j})\Delta\tau_{j}\delta_{t}^{\alpha_{j}}u_{h}^{n},v\bigg)_{D^{k}}
-\varepsilon(p_{h}^{n},v\big)_{D^{k}}
-\big(f(u_{h}^{n}),\frac{\partial}{\partial x}v\big)_{D^{k}}
+\big((\widehat{f}(u_{h}^{n})v^{-})_{k+\frac{1}{2}}-(\widehat{f}(u_{h}^{n})v^{+})_{k-\frac{1}{2}}\big)=0,\\
&\big(p_{h}^{n},\psi\big)_{D^{k}}=\big(\Delta_{(\beta-2)/2}q^{n}_{h},\psi\big)_{D^{k}},\\
&\big(q_{h}^{n},\phi\big)_{D^{k}}=-\big(r_{h}^{n},\frac{\partial \phi}{\partial x}\big)_{D^{k}}+\big((\widehat{r}^{n}_{h}\phi^{-})_{k+\frac{1}{2}}
-(\widehat{r}^{n}_{h}\phi^{+})_{k-\frac{1}{2}}\big),\\
&\big(r_{h}^{n},\eta\big)_{D^{k}}=-\big(u_{h}^{n},\frac{\partial \eta}{\partial x}\big)_{D^{k}}+\big((\widehat{u}^{n}_{h}\eta^{-})_{k+\frac{1}{2}}
-(\widehat{u}^{n}_{h}\eta^{+})_{k-\frac{1}{2}}\big).\\
\end{split}
\end{equation}
The 'hat' terms in the scheme are the so-called numerical fluxes. In order to ensure the stability,
these terms are taken as
\begin{equation}\label{flz}
\begin{split}
\widehat{u}_{h}^{n}=(u_{h}^{n})^{-},\quad\widehat{r}_{h}^{n}=(r_{h}^{n})^{+},\quad \widehat{f}_{h}=\widehat{f}((u_{h}^{n})^{-},(u_{h}^{n})^{+}).
\end{split}
\end{equation}
Note that we can also choose
\begin{equation}\label{1a}
\begin{split}
\widehat{u}_{h}^{n}=(u_{h}^{n})^{+},\quad\widehat{r}_{h}^{n}=(r_{h}^{n})^{-},\quad \widehat{f}_{h}=\widehat{f}((u_{h}^{n})^{-},(u_{h}^{n})^{+}).
\end{split}
\end{equation}

\section{ Stability and error estimates}\label{s3}
 In the following we discuss stability and accuracy of the proposed scheme, for  time distributed order and space-fractional  convection-diffusion problem.
\subsection{ The analysis of stability for fully discrete scheme }
\begin{thm}\label{tt4}
The fully-discrete LDG scheme \eqref{1dv} is  stable, and
\begin{equation}\label{1dvzv}
\begin{split}
\|u_{h}^{n}\|_{L^{2}(\Omega)}\leq C\|u_{h}^{0}\|_{L^{2}(\Omega)},\quad n=1,2,...,M.
\end{split}
\end{equation}
\end{thm}
\textbf{Proof.}  Set $(v,\psi,\phi,\eta)=(u_{h}^{n},p_{h}^{n}-q_{h}^{n}, u_{h}^{n},r_{h}^{n})$ in \eqref{1dv}, and define $\theta(u_{h}^{n})=\int^{u_{h}^{n}}f(s_{h}^{n})ds_{h}^{n}$. Then the following result holds:
\begin{equation}\label{1dv1cv}
\begin{split}
&\bigg(\sum_{j=1}^{S}W(\alpha_{j})\Delta\tau_{j}\delta_{t}^{\alpha_{j}}u_{h}^{n},u_{h}^{n}\bigg)_{D^{k}}
-\varepsilon\big(p_{h}^{n}, u_{h}^{n}\big)_{D^{k}}-(\theta(u_{h}^{n}))^{-}_{k+\frac{1}{2}}
+(\theta(u_{h}^{n}))^{+}_{k-\frac{1}{2}}
+\big((\widehat{f}(u_{h}^{n})u^{-})_{k+\frac{1}{2}}-(\widehat{f}(u_{h}^{n})
(u_{h}^{n})^{+})_{k-\frac{1}{2}}\big)\\
&\quad-\big(p_{h}^{n},q_{h}^{n}\big)_{D^{k}}
+\big(p_{h}^{n},p_{h}^{n}\big)_{D^{k}}+\big(\Delta_{(\beta-2)/2}q^{n}_{h},q_{h}^{n}\big)_{D^{k}}
-\big(\Delta_{(\beta-2)/2}q^{n}_{h},p_{h}^{n}\big)_{D^{k}}+\big(r_{h}^{n},r_{h}^{n}\big)_{D^{k}}
+\big(u_{h}^{n},\frac{\partial r_{h}^{n}}{\partial x}\big)_{D^{k}}\\
&\quad+\big(q_{h}^{n},u_{h}^{n}\big)_{D^{k}}
+\big(r_{h}^{n},\frac{\partial u_{h}^{n}}{\partial x}\big)_{D^{k}}
-\big((\widehat{u}^{n}_{h}(r_{h}^{n})^{-})_{k+\frac{1}{2}}
-(\widehat{u}^{n}_{h}(r_{h}^{n})^{+})_{k-\frac{1}{2}}\big)-\big((\widehat{r}^{n}_{h}(u_{h}^{n})^{-})_{k+\frac{1}{2}}
-(\widehat{r}^{n}_{h}(u_{h}^{n})^{+})_{k-\frac{1}{2}}\big)=0.\\
\end{split}
\end{equation}
Summing over $k$, with the definition \eqref{flz} of the numerical fluxes and with simple algebraic manipulations and, we easily obtain
\begin{equation}\label{1dv1}
\begin{split}
&\bigg(\sum_{j=1}^{S}W(\alpha_{j})\Delta\tau_{j}\delta_{t}^{\alpha_{j}}u_{h}^{n},u_{h}^{n}\bigg)
-\varepsilon\big(p_{h}^{n}, u_{h}^{n}\big)-\sum_{k=1}^{K}\big((\theta(u_{h}^{n}))^{-}_{k+\frac{1}{2}}
-(\theta(u_{h}^{n}))^{+}_{k-\frac{1}{2}}\big)
+\sum_{k=1}^{K}\big((\widehat{f}(u_{h}^{n})u^{-})_{k+\frac{1}{2}}-(\widehat{f}(u_{h}^{n})
(u_{h}^{n})^{+})_{k-\frac{1}{2}}\big)\\
&\quad-\big(p_{h}^{n},q_{h}^{n}\big)
+\big(p_{h}^{n},p_{h}^{n}\big)+\big(\Delta_{(\beta-2)/2}q^{n}_{h},q_{h}^{n}\big)
-\big(\Delta_{(\beta-2)/2}q^{n}_{h},p_{h}^{n}\big)
+\big(r_{h}^{n},r_{h}^{n}\big)+\big(q_{h}^{n},u_{h}^{n}\big)
=0.\\
\end{split}
\end{equation}
From the properties of the monotone flux, we know that $\widehat{f}((u_{h}^{n})^{-},(u_{h}^{n})^{+})$ nondecreasing function of its first argument and a nonincreasing function of its second argument. Hence, we have
\begin{equation}\label{1dv1}
\begin{split}
\sum_{k=1}^{K}\big((\widehat{f}(u_{h}^{n})u^{-})_{k+\frac{1}{2}}
-(\widehat{f}(u_{h}^{n})(u_{h}^{n})^{+})_{k-\frac{1}{2}}\big)-
\sum_{k=1}^{K}\big((\theta(u_{h}^{n}))^{-}_{k+\frac{1}{2}}-(\theta(u_{h}^{n}))^{+}_{k-\frac{1}{2}}\big)>0.
\end{split}
\end{equation}
This implies that
\begin{equation}\label{1dv1}
\begin{split}
&\bigg(\sum_{j=1}^{S}W(\alpha_{j})\Delta\tau_{j}\delta_{t}^{\alpha_{j}}u_{h}^{n},u_{h}^{n}\bigg)
+\big(r_{h}^{n},r_{h}^{n}\big)
+\big(p_{h}^{n},p_{h}^{n}\big)+\big(\Delta_{(\beta-2)/2}q^{n}_{h},q_{h}^{n}\big)\\
&\leq\big(\Delta_{(\beta-2)/2}q^{n}_{h},p_{h}^{n}\big)
-\big(q_{h}^{n},u_{h}^{n}\big)+\big(p_{h}^{n},q_{h}^{n}\big)
+\varepsilon\big(p_{h}^{n}, u_{h}^{n}\big).\\
\end{split}
\end{equation}
Employing Young's inequality and  Lemma \ref{lga2}, we obtain
\begin{equation}\label{1dv1}
\begin{split}
&\bigg(\sum_{j=1}^{S}W(\alpha_{j})\Delta\tau_{j}\delta_{t}^{\alpha_{j}}u_{h}^{n},u_{h}^{n}\bigg)
+\|r_{h}^{n}\|^{2}_{L^{2}(\Omega)}+\big(\Delta_{(\beta-2)/2}q^{n}_{h},q_{h}^{n}\big)\leq c\|u_{h}^{n}\|^{2}_{L^{2}(\Omega)}+c_{1}\|q_{h}^{n}\|^{2}_{L^{2}(\Omega)}.\\
\end{split}
\end{equation}
Recalling Lemma \ref{lga2}, we obtain
\begin{equation}\label{1dv1}
\begin{split}
&\bigg(\sum_{j=1}^{S}W(\alpha_{j})\Delta\tau_{j}\delta_{t}^{\alpha_{j}}u_{h}^{n},u_{h}^{n}\bigg)+\|r_{h}^{n}\|^{2}_{L^{2}(\Omega)}\leq c\|u_{h}^{n}\|^{2}_{L^{2}(\Omega)}.\\
\end{split}
\end{equation}
It then follows that
\begin{equation}\label{1dv1}
\begin{split}
\bigg(\sum_{j=1}^{S}
\frac{W(\alpha_{j})\Delta\tau_{j}}{\lambda_{j}}u^{n},u^{n}_{h}\bigg)\leq&
\bigg(\sum_{j=1}^{S}
\frac{W(\alpha_{j})\Delta\tau_{j}}{\lambda_{j}}\sum_{l=1}^{n-1}
(a_{n-l-1}^{\alpha_{j}}-a_{n-l}^{\alpha_{j}})u^{l}_{h},u^{n}_{h}\bigg)
\\
&\quad+
\bigg(\sum_{j=1}^{S}
\frac{W(\alpha_{j})\Delta\tau_{j}}{\lambda_{j}}a_{n-1}^{\alpha_{j}}u^{0}_{h},u^{n}_{h}\bigg)+ c\|u_{h}^{n}\|^{2}_{L^{2}(\Omega)}.\\
\end{split}
\end{equation}
Using Cauchy-Schwarz inequality, we obtain
\begin{equation}\label{1dv1}
\begin{split}
\|u_{h}^{n}\|^{2}_{L^{2}(\Omega)}\leq &c_{1}\sum_{j=1}^{S}
\frac{W(\alpha_{j})\Delta\tau_{j}}{\lambda_{j}}Q\sum_{l=1}^{n-1}
(a_{n-l-1}^{\alpha_{j}}-a_{n-l}^{\alpha_{j}})\|u_{h}^{l}\|_{L^{2}(\Omega)}\|u_{h}^{n}\|_{L^{2}(\Omega)}\\
&\quad+c_{2} \sum_{j=1}^{S}
\frac{W(\alpha_{j})\Delta\tau_{j}}{\lambda_{j}}Qa_{n-1}^{\alpha_{j}}\|u_{h}^{0}\|_{L^{2}(\Omega)}\|u_{h}^{n}\|_{L^{2}(\Omega)}
+cQ\|u_{h}^{n}\|^{2}_{L^{2}(\Omega)},\\
\end{split}
\end{equation}
where $Q=\frac{1}{\sum_{j=1}^{S}\frac{W(\alpha_{j})\Delta\tau_{j}}{\lambda_{j}}}$ and provided $c$ is sufficiently small such that $1-cQ>0$, we obtain that

\begin{equation}\label{1dv1zx}
\begin{split}
&\|u_{h}^{n}\|_{L^{2}(\Omega)}\leq C\bigg(\sum_{j=1}^{S}
\frac{W(\alpha_{j})\Delta\tau_{j}}{\lambda_{j}}Q\sum_{l=1}^{n-1}
(a_{n-l-1}^{\alpha_{j}}-a_{n-l}^{\alpha_{j}})\|u_{h}^{l}\|_{L^{2}(\Omega)}+ \sum_{j=1}^{S}
\frac{W(\alpha_{j})\Delta\tau_{j}}{\lambda_{j}}Qa_{n-1}^{\alpha_{j}}\|u_{h}^{0}\|_{L^{2}(\Omega)}\bigg).\\
\end{split}
\end{equation}
Obviously the theorem holds for $n = 0$. Assume that it is valid for $n = 1, 2, ..., m-1$. Then,
by \eqref{1dv1zx}, we have
\begin{equation}\label{1dv1}
\begin{split}
\|u_{h}^{m}\|_{L^{2}(\Omega)}&\leq C\bigg(\sum_{j=1}^{S}
\frac{W(\alpha_{j})\Delta\tau_{j}}{\lambda_{j}}Q\sum_{l=1}^{m-1}
(a_{n-l-1}^{\alpha_{j}}-a_{n-l}^{\alpha_{j}})\|u_{h}^{l}\|_{L^{2}(\Omega)}+ \sum_{j=1}^{S}
\frac{W(\alpha_{j})\Delta\tau_{j}}{\lambda_{j}}Qa_{n-1}^{\alpha_{j}}\|u_{h}^{0}\|_{L^{2}(\Omega)}\bigg)\\
&\leq C\bigg(\sum_{j=1}^{S}
\frac{W(\alpha_{j})\Delta\tau_{j}}{\lambda_{j}}Q\sum_{l=1}^{m-1}
(a_{n-l-1}^{\alpha_{j}}-a_{n-l}^{\alpha_{j}})\|u_{h}^{0}\|_{L^{2}(\Omega)}+ \sum_{j=1}^{S}
\frac{W(\alpha_{j})\Delta\tau_{j}}{\lambda_{j}}Qa_{n-1}^{\alpha_{j}}\|u_{h}^{0}\|_{L^{2}(\Omega)}\bigg)\\
&=C\|u_{h}^{0}\|_{L^{2}(\Omega)}. \quad\Box\\
\end{split}
\end{equation}

\subsection{Error estimates}
In order to obtain the error estimate to smooth solutions
for the considered fully discrete LDG scheme \eqref{1dv}, we need to first obtain the error equation.\\
It is easy to verify that the exact solution of \eqref{25n} satisfies
\begin{equation}\label{1dv1ffx}
\begin{split}
&\bigg(\sum_{j=1}^{S}W(\alpha_{j})\Delta\tau_{j}\delta_{t}^{\alpha_{j}}u^{n},v\bigg)_{D^{k}}
-\varepsilon\big(p^{n},v\big)_{D^{k}}
+\big(\gamma(x)^{n},v\big)_{D^{k}}-\big(f(u^{n}),\frac{\partial}{\partial x}v\big)_{D^{k}} \\
&\quad\quad\quad+\big((\widehat{f}(u^{n})v^{-})_{k+\frac{1}{2}}
-(\widehat{f}(u^{n})v^{+})_{k-\frac{1}{2}}\big)=0,\\
&\big(p^{n},\psi\big)_{D^{k}} =\big(\Delta_{(\beta-2)/2}q^{n},\psi\big)_{D^{k}} ,\\
&\big(q^{n},\phi\big)_{D^{k}} =-\big(r^{n},\frac{\partial \phi}{\partial x}\big)_{D^{k}} +\big((\widehat{r}^{n}\phi^{-})_{k+\frac{1}{2}}
-(\widehat{r}^{n}\phi^{+})_{k-\frac{1}{2}}\big),\\
&\big(r^{n},\eta\big)_{D^{k}} =-\big(u^{n},\frac{\partial \eta}{\partial x}\big)_{D^{k}} +\big((\widehat{u}^{n}\eta^{-})_{k+\frac{1}{2}}
-(\widehat{u}^{n}\eta^{+})_{k-\frac{1}{2}}\big).\\
\end{split}
\end{equation}
Subtracting equation \eqref{1dv} from \eqref{1dv1ffx}, we can obtain the error equation
\begin{equation}\label{1dv1hj}
\begin{split}
\bigg(\sum_{j=1}^{S}&W(\alpha_{j})\Delta\tau_{j}\delta_{t}^{\alpha_{j}}(u^{n}-u^{n}_{h}),v\bigg)_{D^{k}}
-\varepsilon\big(p^{n}-p^{n}_{h},v\big)_{D^{k}}
+\big(\gamma(x)^{n},v\big)_{D^{k}} -\big(f(u^{n})-f(u^{n}_{h}),\frac{\partial}{\partial x}v\big)_{D^{k}}
\\
&+\big(((\widehat{f}(u^{n})-\widehat{f}(u^{n}_{h}))v^{-})_{k+\frac{1}{2}}-((\widehat{f}(u^{n})-
\widehat{f}(u^{n}_{h}))v^{+})_{k-\frac{1}{2}}\big)
+\big(p^{n}-p^{n}_{h},\psi\big)_{D^{k}} +\big(q^{n}-q^{n}_{h},\phi\big)_{D^{k}}
\\
& -\big(\Delta_{(\beta-2)/2}(q^{n}-q^{n}_{h}),\psi\big)_{D^{k}}+\big(r^{n}-r^{n}_{h},\frac{\partial \phi}{\partial x}\big)_{D^{k}} -((\widehat{r}^{n}-\widehat{r}^{n}_{h})\phi^{-})_{k+\frac{1}{2}}
+((\widehat{r}^{n}-\widehat{r}^{n}_{h})\phi^{+})_{k-\frac{1}{2}}
 \\
&+\big(r^{n}-r^{n}_{h},\eta\big)_{D^{k}}+\big(u^{n}-u^{n}_{h},\frac{\partial \eta}{\partial x}\big)_{D^{k}} -((\widehat{u}^{n}-\widehat{u}^{n}_{h})\eta^{-})_{k+\frac{1}{2}}
+((\widehat{u}^{n}-\widehat{u}^{n}_{h})\eta^{+})_{k-\frac{1}{2}}=0,\\
\end{split}
\end{equation}
where
\begin{equation}\label{1dv1ffasnm}
\begin{split}
|\gamma(x)^{n}|=|\mathcal{O}((\Delta t)^{2-\alpha_{j}}+\theta^{2})|\leq c((\Delta t)^{1+\frac{\theta}{2}}+\theta^{2}),
\end{split}
\end{equation}
 such that
\begin{equation}\label{1dv1ffaszz}
\begin{split}
1+\frac{\theta}{2}=2-S\theta+\frac{\theta}{2}\leq 2-\alpha_{j}=2-j\theta+\frac{\theta}{2}\leq2-\theta+\frac{\theta}{2}=2-\frac{\theta}{2}.
\end{split}
\end{equation}
For the error estimate, we define special projections, $\mathcal{P}$ and $\mathcal{P}^{\pm}$ into $V_{h}^{k}$.  For
all the elements, $D^{k}$, $k = 1, 2, ... , K$ are defined to satisfy
\begin{equation}\label{prh}
\begin{split}
&(\mathcal{P}u-u,v)_{D^{k}}=0,\quad\forall  v\in\mathbb{P}_{N}^{k}(D^{k}),\\
&(\mathcal{P}^{\pm}u-u,v)_{D^{k}}=0,\quad\forall  v\in\mathbb{P}_{N}^{k-1}(D^{k}),\quad \mathcal{P}^{\pm}u_{k+\frac{1}{2}}=u(x_{k+\frac{1}{2}}^{\pm}).\\
\end{split}
\end{equation}
Denoting
\begin{equation}\label{91hgh}
\begin{split}
&\pi^{n}=\mathcal{P}^{-}u^{n}-u_{h}^{n},\quad \pi^{e}_{n}=\mathcal{P}^{-}u^{n}-u^{n},\quad \sigma^{n}=\mathcal{P}p^{n}-p_{h}^{n},\quad \sigma^{e}_{n}=\mathcal{P}p^{n}-p^{n},\quad\varphi^{n}=\mathcal{P}q^{n}-q_{h}^{n},\\ &\varphi^{e}_{n}=\mathcal{P}q^{n}-q^{n},\quad\psi^{n}=\mathcal{P}^{+}r^{n}-r_{h}^{n},\quad \psi^{e}_{n}=\mathcal{P}^{+}r^{n}-r^{n}.
\end{split}
\end{equation}
For the special projections mentioned above, we have, by the standard approximation
theory \cite{Ciarlet:2002:FEM:581834}, that
\begin{equation}\label{sth}
\begin{split}
&\|\pi^{e}\|_{L^{2}(\Omega)}+h\|\pi^{e}\|_{\infty}+h^{\frac{1}{2}}\| \pi\|_{\Gamma_{h}}\leq  Ch^{N+1}.\\
\end{split}
\end{equation}
where $\pi^{e}=\mathcal{P}^{\pm}u^{n}-u_{h}^{n}$ or $\pi^{e}=\mathcal{P}u^{n}-u_{h}^{n}$. The positive constant $C$, solely depending on $u^{n}$, is independent of $h$. $\Gamma_{h}$ denotes the set of boundary points of all elements $D^{k}$.
\begin{thm}\label{th1} (Diffusion without convection $f(u)=0$).
 Let $u(x,t_{n})$ be the exact solution of the problem \eqref{25n}, which
is sufficiently smooth with bounded derivatives, let $u_{h}^{n}$ be the numerical solution
of the fully discrete LDG scheme \eqref{1dv}, then there holds the following error
estimates:
\begin{equation}\label{tt7fz}
\begin{split}
&\|u(x,t_{n})-u_{h}^{n}\|_{L^{2}(\Omega)}\leq C(h^{N+1}+(\Delta t)^{1+\frac{\theta}{2}}+\theta^{2}).\\
\end{split}
\end{equation}
\end{thm}
\textbf{Proof}. From \eqref{1dv1hj}, we can obtain the error equation
\begin{equation}\label{1dv1gg}
\begin{split}
\bigg(\sum_{j=1}^{S}&W(\alpha_{j})\Delta\tau_{j}\delta_{t}^{\alpha_{j}}(u^{n}-u^{n}_{h}),v\bigg)_{D^{k}}
-\varepsilon\big(p^{n}-p^{n}_{h},v\big)_{D^{k}}
+\big(\gamma(x)^{n},v\big)_{D^{k}} \\
&
+\big(p^{n}-p^{n}_{h},\psi\big)_{D^{k}} -\big(\Delta_{(\beta-2)/2}(q^{n}-q^{n}_{h}),\psi\big)_{D^{k}}
+\big(q^{n}-q^{n}_{h},\phi\big)_{D^{k}} +\big(r^{n}-r^{n}_{h},\frac{\partial \phi}{\partial x}\big)_{D^{k}} \\
&-((\widehat{r}^{n}-\widehat{r}^{n}_{h})\phi^{-})_{k+\frac{1}{2}}
+((\widehat{r}^{n}-\widehat{r}^{n}_{h})\phi^{+})_{k-\frac{1}{2}}
+\big(r^{n}-r^{n}_{h},\eta\big)_{D^{k}} +\big(u^{n}-u^{n}_{h},\frac{\partial \eta}{\partial x}\big)_{D^{k}} \\
&-((\widehat{u}^{n}-\widehat{u}^{n}_{h})\eta^{-})_{k+\frac{1}{2}}
+((\widehat{u}^{n}-\widehat{u}^{n}_{h})\eta^{+})_{k-\frac{1}{2}}=0.\\
\end{split}
\end{equation}
Using  \eqref{91hgh}, the error equation \eqref{1dv1gg} can be written
\begin{equation}\label{1dv1}
\begin{split}
\bigg(\sum_{j=1}^{S}&W(\alpha_{j})\Delta\tau_{j}\delta_{t}^{\alpha_{j}}(\pi^{n}-\pi^{e}_{n}),v\bigg)_{D^{k}}
-\varepsilon\big(\sigma^{n}-\sigma^{e}_{n},v\big)_{D^{k}}
+\big(\gamma(x)^{n},v\big)_{D^{k}} \\
&
+\big(\sigma^{n}-\sigma^{e}_{n},\psi\big)_{D^{k}} -
\big(\Delta_{(\beta-2)/2}(\varphi^{n}-\varphi^{e}_{n}),\psi\big)_{D^{k}}
+\big(\varphi^{n}-\varphi^{e}_{n},\phi\big)_{D^{k}} +\big(\psi^{n}-\psi^{e}_{n},\frac{\partial \phi}{\partial x}\big)_{D^{k}} \\
&-((\psi^{n}-\psi^{e}_{n})^{-}\phi^{-})_{k+\frac{1}{2}}
+((\psi^{n}-\psi^{e}_{n})^{-}\phi^{+})_{k-\frac{1}{2}}
+\big(\psi^{n}-\psi^{e}_{n},\eta\big)_{D^{k}} +\big(\pi^{n}-\pi^{e}_{n},\frac{\partial \eta}{\partial x}\big)_{D^{k}} \\
&-((\pi^{n}-\pi^{e}_{n})^{+}\eta^{-})_{k+\frac{1}{2}}
+((\pi^{n}-\pi^{e}_{n})^{+}\eta^{+})_{k-\frac{1}{2}}=0,\\
\end{split}
\end{equation}
and taking the test functions
\begin{equation}\label{91h}
\begin{split}
v=\pi^{n},\quad  \psi=\sigma^{n}-\varphi^{n},\quad\phi=\pi^{n},\quad\eta=\psi^{n},
\end{split}
\end{equation}
 we obtain
\begin{equation}\label{1dv1ggh}
\begin{split}
\bigg(\sum_{j=1}^{S}&W(\alpha_{j})\Delta\tau_{j}\delta_{t}^{\alpha_{j}}(\pi^{n}-\pi^{e}_{n}),\pi^{n}\bigg)
-\varepsilon\big(\sigma^{n}-\sigma^{e}_{n},\pi^{n}\big)
+\big(\gamma(x)^{n},\pi^{n}\big)\\
&+\big(\sigma^{n}-\sigma^{e}_{n},-\varphi^{n}+\sigma^{n}\big)-
\big(\Delta_{(\beta-2)/2}(\varphi^{n}-\varphi^{e}_{n}),-\varphi^{n}+\sigma^{n}\big)
+\big(\varphi^{n}-\varphi^{e}_{n},\pi^{n}\big)\\
&-\big(\psi^{e}_{n},\frac{\partial \pi^{n}}{\partial x}\big)+\sum_{k=1}^{K}\big(((\psi^{e}_{n})^{-}(\pi^{n})^{-})_{k+\frac{1}{2}}
-((\psi^{e}_{n})^{-}(\pi^{n})^{+})_{k-\frac{1}{2}}\big)
+\big(\psi^{n}-\psi^{e}_{n},\psi^{n}\big)\\
&-\big(\pi^{e}_{n},\frac{\partial \psi^{n}}{\partial x}\big)+\sum_{k=1}^{K}\big(((\pi^{e}_{n})^{+}(\psi^{n})^{-})_{k+\frac{1}{2}}
-((\pi^{e}_{n})^{+}(\psi^{n})^{+})_{k-\frac{1}{2}}\big)=0,\\
\end{split}
\end{equation}
by the properties of the projection $P^{+}$ and $P^{-}$ we obtain
\begin{equation}\label{1dv1ff}
\begin{split}
\bigg(\sum_{j=1}^{S}&W(\alpha_{j})\Delta\tau_{j}\delta_{t}^{\alpha_{j}}(\pi^{n}-\pi^{e}_{n}),\pi^{n}\bigg)
-\varepsilon\big(\sigma^{n}-\sigma^{e}_{n},\pi^{n}\big)
+\big(\gamma(x)^{n},\pi^{n}\big)\\
&
+\big(\sigma^{n}-\sigma^{e}_{n},-\varphi^{n}+\sigma^{n}\big)-
\big(\Delta_{(\beta-2)/2}(\varphi^{n}-\varphi^{e}_{n}),-\varphi^{n}+\sigma^{n}\big)
+\big(\varphi^{n}-\varphi^{e}_{n},\pi^{n}\big)\\
&+\sum_{k=1}^{K}\big(((\psi^{e}_{n})^{-}(\pi^{n})^{-})_{k+\frac{1}{2}}
-((\psi^{e}_{n})^{-}(\pi^{n})^{+})_{k-\frac{1}{2}}\big)
+\big(\psi^{n}-\psi^{e}_{n},\psi^{n}\big)\\
&+\sum_{k=1}^{K}\big(((\pi^{e}_{n})^{+}(\psi^{n})^{-})_{k+\frac{1}{2}}
-((\pi^{e}_{n})^{+}(\psi^{n})^{+})_{k-\frac{1}{2}}\big)=0.\\
\end{split}
\end{equation}
Employing Young's inequality and Lemma \ref{lga2} and  the interpolation property \eqref{sth} and \eqref{1dv1ffasnm}, we obtain
\begin{equation}\label{1dv1}
\begin{split}
\bigg(\sum_{j=1}^{S}&W(\alpha_{j})\Delta\tau_{j}\delta_{t}^{\alpha_{j}}\pi^{n},\pi^{n}\bigg)
+\big(\sigma^{n},\sigma^{n}\big)+
\big(\Delta_{(\beta-2)/2}\varphi^{n},\varphi^{n}\big)
+\big(\psi^{n},\psi^{n}\big)\\
&\leq
C(h^{2N+2}+(\Delta t)^{4+\theta}+\theta^{4})+\bigg(\sum_{j=1}^{S}W(\alpha_{j})\Delta\tau_{j}\delta_{t}^{\alpha_{j}}\pi^{e}_{n}
,\pi^{n}\bigg)+c_{2}\|\sigma^{n}\|^{2}_{L^{2}(\Omega)}+ c_{3}\|\varphi^{n}\|^{2}_{L^{2}(\Omega)}\\
&\quad+c\|\pi^{n}\|^{2}_{L^{2}(\Omega)}
+c_{1}\|\psi^{n}\|^{2}_{L^{2}(\Omega)},
\end{split}
\end{equation}
by using Lemma \ref{lga2z}, \eqref{91hgh} and the interpolation property \eqref{sth}, we get
\begin{equation}\label{sthzz}
\begin{split}
&\|\delta_{t}^{\alpha}(\mathcal{P}^{+}u(x,t_{n})-u(x,t_{n}))\|_{L^{2}(\Omega)}\leq  C\big(h^{N+1}+(\Delta t)^{2-\alpha}\big).\\
\end{split}
\end{equation}
From \eqref{7zz}, \eqref{1dv1ffaszz} and \eqref{sthzz}, we obtain
\begin{equation}\label{sthzzc}
\begin{split}
&\Biggl\|\sum_{j=1}^{S}W(\alpha_{j})\Delta\tau_{j}\delta_{t}^{\alpha_{j}}(\mathcal{P}^{+}u(x,t_{n})-u(x,t_{n}))
\Biggr\|_{L^{2}(\Omega)}\leq  C\big(h^{N+1}+(\Delta t)^{1+\frac{\theta}{2}}+\theta^{2}\big).\\
\end{split}
\end{equation}
Hence
\begin{equation}\label{1dv1}
\begin{split}
&\bigg(\sum_{j=1}^{S}W(\alpha_{j})\Delta\tau_{j}\delta_{t}^{\alpha_{j}}\pi^{n},\pi^{n}\bigg)
+\big(\sigma^{n},\sigma^{n}\big)+
\big(\Delta_{(\beta-2)/2}\varphi^{n},\varphi^{n}\big)
+\big(\psi^{n},\psi^{n}\big)\\
&\leq
C\big(h^{2N+2}+(\Delta t)^{2+\theta}+\theta^{4}\big)+c_{2}\|\sigma^{n}\|^{2}_{L^{2}(\Omega)}+ c_{3}\|\varphi^{n}\|^{2}_{L^{2}(\Omega)}+c\|\pi^{n}\|^{2}_{L^{2}(\Omega)}
+c_{1}\|\psi^{n}\|^{2}_{L^{2}(\Omega)}.
\end{split}
\end{equation}
Recalling Lemma \ref{lg} and provided $c_{i},\,\,i=1,2$ are sufficiently small such that $c_{i}\leq1$, we obtain
\begin{equation}\label{1dv1}
\begin{split}
&\bigg(\sum_{j=1}^{S}W(\alpha_{j})\Delta\tau_{j}\delta_{t}^{\alpha_{j}}\pi^{n},\pi^{n}\bigg)
\leq
C\big(h^{2N+2}+(\Delta t)^{2+\theta}+\theta^{4}\big)+c\|\pi^{n}\|^{2}_{L^{2}(\Omega)}.
\end{split}
\end{equation}
It then follows that
\begin{equation}\label{1dv1}
\begin{split}
\bigg(\sum_{j=1}^{S}
\frac{W(\alpha_{j})\Delta\tau_{j}}{\lambda_{j}}\pi^{n},\pi^{n}_{h}\bigg)\leq&
\bigg(\sum_{j=1}^{S}
\frac{W(\alpha_{j})\Delta\tau_{j}}{\lambda_{j}}\sum_{l=1}^{n-1}
(a_{n-l-1}^{\alpha_{j}}-a_{n-l}^{\alpha_{j}})\pi^{l},\pi^{n}\bigg)\\
&+
\bigg(\sum_{j=1}^{S}
\frac{W(\alpha_{j})\Delta\tau_{j}}{\lambda_{j}}a_{n-1}^{\alpha_{j}}\pi^{0},\pi^{n}\bigg)+ c\|\pi^{n}\|^{2}_{L^{2}(\Omega)}+C\big(h^{2N+2}+(\Delta t)^{2+\theta}+\theta^{4}\big).\\
\end{split}
\end{equation}
Employing Young's inequality, we obtain
\begin{equation}\label{1dv1}
\begin{split}
\|\pi^{n}\|^{2}_{L^{2}(\Omega)}\leq &\sum_{j=1}^{S}
\frac{W(\alpha_{j})\Delta\tau_{j}}{\lambda_{j}}Q\sum_{l=1}^{n-1}
(a_{n-l-1}^{\alpha_{j}}-a_{n-l}^{\alpha_{j}})\|\pi^{l}\|^{2}_{L^{2}(\Omega)}
+\frac{1}{4}\sum_{j=1}^{S}
\frac{W(\alpha_{j})\Delta\tau_{j}}{\lambda_{j}}Q(a_{0}^{\alpha_{j}}-a_{n-1}^{\alpha_{j}})
\|\pi^{n}\|^{2}_{L^{2}(\Omega)}\\
&+ \sum_{j=1}^{S}
\frac{W(\alpha_{j})\Delta\tau_{j}}{\lambda_{j}}Qa_{n-1}^{\alpha_{j}}\|\pi^{0}\|^{2}_{L^{2}(\Omega)}
+\sum_{j=1}^{S}
\frac{W(\alpha_{j})\Delta\tau_{j}}{4\lambda_{j}}Qa_{n-1}^{\alpha_{j}}\|\pi^{n}\|^{2}_{L^{2}(\Omega)}\\
&+cQ\|\pi^{n}\|^{2}_{L^{2}(\Omega)}
+CQ\big(h^{2N+2}+(\Delta t)^{2+\theta}+\theta^{4}\big).\\
\end{split}
\end{equation}
Notice the facts that
\begin{equation}\label{1dvc42}
\begin{split}
\|\pi^{0}\|_{L^{2}(\Omega)}\leq Ch^{N+1}.
\end{split}
\end{equation}
Thus,
\begin{equation}\label{1dv1}
\begin{split}
\|\pi^{n}\|^{2}_{L^{2}(\Omega)}\leq &\sum_{j=1}^{S}
\frac{W(\alpha_{j})\Delta\tau_{j}}{\lambda_{j}}Q\sum_{l=1}^{n-1}
(a_{n-l-1}^{\alpha_{j}}-a_{n-l}^{\alpha_{j}})\|\pi^{l}\|^{2}_{L^{2}(\Omega)}
+(cQ+\frac{1}{4})\sum_{j=1}^{S}
\frac{W(\alpha_{j})\Delta\tau_{j}}{\lambda_{j}}Q\|\pi^{n}\|^{2}_{L^{2}(\Omega)}\\
&+C\sum_{j=1}^{S}
\frac{W(\alpha_{j})\Delta\tau_{j}}{\lambda_{j}}Qa_{n-1}^{\alpha_{j}}h^{2N+2}+C\sum_{j=1}^{S}
\frac{W(\alpha_{j})\Delta\tau_{j}}{\lambda_{j}}Qa_{n-1}^{\alpha_{j}}\big(h^{2N+2}+(\Delta t)^{2+\theta}+\theta^{4}\big),\\
\end{split}
\end{equation}
  provided $c$ is sufficiently small such that $\frac{3}{4}-cQ>0$, we obtain that
\begin{equation}\label{1dv1mn}
  \begin{split}
&\|\pi^{n}\|^{2}_{L^{2}(\Omega)}\leq C\bigg(\sum_{j=1}^{S}
\frac{W(\alpha_{j})\Delta\tau_{j}}{\lambda_{j}}Q\sum_{l=1}^{n-1}
(a_{n-l-1}^{\alpha_{j}}-a_{n-l}^{\alpha_{j}})\|\pi^{l}\|^{2}_{L^{2}(\Omega)}
+\sum_{j=1}^{S}
\frac{W(\alpha_{j})\Delta\tau_{j}}{\lambda_{j}}Qa_{n-1}^{\alpha_{j}}\big(h^{2N+2}+(\Delta t)^{2+\theta}+\theta^{4}\big)\bigg).\\
\end{split}
\end{equation}
Obviously the theorem holds for $n = 0$. Assume that it is valid for $n = 1, 2, ..., m-1$. Then,
by \eqref{1dv1mn}, we have
\begin{equation}\label{1dv1}
  \begin{split}
&\|\pi^{m}\|^{2}_{L^{2}(\Omega)}\leq C\bigg(\sum_{j=1}^{S}
\frac{W(\alpha_{j})\Delta\tau_{j}}{\lambda_{j}}Q\sum_{l=1}^{m-1}
(a_{n-l-1}^{\alpha_{j}}-a_{n-l}^{\alpha_{j}})\|\pi^{l}\|^{2}_{L^{2}(\Omega)}
+\sum_{j=1}^{S}
\frac{W(\alpha_{j})\Delta\tau_{j}}{\lambda_{j}}Qa_{n-1}^{\alpha_{j}}\big(h^{2N+2}+(\Delta t)^{2+\theta}+\theta^{4}\big)\bigg)\\
&\leq C\bigg(\sum_{j=1}^{S}
\frac{W(\alpha_{j})\Delta\tau_{j}}{\lambda_{j}}Q\sum_{l=1}^{m-1}
(a_{n-l-1}^{\alpha_{j}}-a_{n-l}^{\alpha_{j}})(h^{2N+2}+(\Delta t)^{2+\theta}+\theta^{4})
+\sum_{j=1}^{S}
\frac{W(\alpha_{j})\Delta\tau_{j}}{\lambda_{j}}Qa_{n-1}^{\alpha_{j}}\big(h^{2N+2}+(\Delta t)^{2+\theta}+\theta^{4}\big)\bigg)
\\
&= C\big(h^{2N+2}+(\Delta t)^{2+\theta}+\theta^{4}\big).
\end{split}
\end{equation}
Finally, by using triangle inequality and standard approximation theory, we get
 \begin{equation}\label{tt7fzx}
\begin{split}
&\|u(x,t_{m})-u_{h}^{m}\|_{L^{2}(\Omega)}\leq C(h^{N+1}+(\Delta t)^{1+\frac{\theta}{2}}+\theta^{2}).\quad\Box\\
\end{split}
\end{equation}
For the more general fractional convection-diffusion problem, we introduce a few
results and then give the error estimate.

\begin{lem}\label{as1}  (see \cite{zhang2004error}). For any piecewise smooth function $\pi \in L^{2}(\Omega)$, on each cell boundary point we define
\begin{equation}\label{91}
\begin{split}
\kappa(\widehat{f};\pi)\equiv \kappa(\widehat{f};\pi^{-},\pi^{+})=\left\{
  \begin{array}{ll}
    [w]^{-1}(f(\pi)-\widehat{f}(\pi)), & \hbox{if $[\pi]\neq0$;} \\
    \frac{1}{2}|f^{'}(\overline{\pi})|, & \hbox{if $[\pi]=0$,}
  \end{array}
\right.
\end{split}
\end{equation}
\end{lem}
where $\widehat{f}(\pi) \equiv \widehat{f}(\pi^{−},\pi^{+})$ is a monotone numerical flux consistent with the given flux
$f$. Then $\kappa(\widehat{f},\pi)$ is nonnegative and bounded for any ($\pi^{−},\pi^{+}) \in \mathbb{R}$. \\
To estimate the nonlinear part,  we can write it into the following
form
\begin{equation}\label{91bm44j}
\begin{split}
\sum_{k=1}^{K}\mathcal{H}_{k}(f;u,u_{h};\pi)=\sum_{k=1}^{K}\big(f(u)-f(u_{h}),\frac{\partial}{\partial x}\pi\big)_{D^{k}}+\sum_{k=1}^{K}((f(u)-f(u_{h})[\pi])_{k+\frac{1}{2}}+\sum_{k=1}^{K}((f(u_{h})
-\widehat{f})[\pi])_{k+\frac{1}{2}}.\\
\end{split}
\end{equation}
We can rewrite \eqref{91bm44j} as:
\begin{equation}\label{91bm44}
\begin{split}
\sum_{k=1}^{K}\mathcal{H}_{k}(f;u,u_{h};\pi)=\sum_{k=1}^{K}\big(f(u)-f(u_{h}),\frac{\partial}{\partial x}\pi\big)_{D^{k}}+\sum_{k=1}^{K}((f(u)-f(\{u_{h}\})[\pi])_{k+\frac{1}{2}}+\sum_{k=1}^{K}((f(\{u_{h}\})
-\widehat{f})[\pi])_{k+\frac{1}{2}}.\\
\end{split}
\end{equation}
\begin{lem}\label{g91g} (see \cite{zhang2004error}) For $\mathcal{H}_{k}(f;u,u_{h};\pi)$ defined above, we have the following
estimate:
\begin{equation}\label{91gf}
\begin{split}
\sum_{k=1}^{K}\mathcal{H}_{k}(f;u,u_{h};v)\leq&-\frac{1}{4}\kappa(\widehat{f};u_{h})
+(C+C_{*}(\|v\|_{\infty}+h^{-1}\|e_{u}\|_{\infty}^{2}))\|v\|^{2}+(C+C_{*}h^{-1}\|e_{u}\|_{\infty}^{2})h^{2N+1}.
\end{split}
\end{equation}
\end{lem}
To deal with the nonlinearity of the flux $f(u)$, we make the following assumption
for $h$ small enough and $k \geq 1$, which can be verified \cite{xu2007error}:
\begin{equation}\label{91}
\begin{split}
\|e_{u}\|=\|u-u_{h}\|\leq h.
\end{split}
\end{equation}
\begin{thm}\label{th2}
 Let $u(x,t_{n})$ be the exact solution of the problem \eqref{25n}, which
is sufficiently smooth with bounded derivatives, let $u_{h}^{n}$ be the numerical solution
of the fully discrete LDG scheme \eqref{1dv}, then there holds the following error
estimates:
\begin{equation}\label{tt7mk}
\begin{split}
&\|u(x,t_{n})-u_{h}^{n}\|_{L^{2}(\Omega)}\leq C\big(h^{N+\frac{1}{2}}+(\Delta t)^{1+\frac{\theta}{2}}+\theta^{2}\big).\\
\end{split}
\end{equation}
\end{thm}
\textbf{Proof}. Using  \eqref{91hgh}, the error equation \eqref{1dv1hj} can be written
\begin{equation}\label{1dv1hj2}
\begin{split}
\bigg(\sum_{j=1}^{S}&W(\alpha_{j})\Delta\tau_{j}\delta_{t}^{\alpha_{j}}(\pi^{n}-\pi^{e}_{n}),v\bigg)
-\varepsilon\big(\sigma^{n}-\sigma^{e}_{n},v\big)
+\big(\gamma(x)^{n},v\big)-\sum_{k=1}^{K}\mathcal{H}_{k}(f;u,u_{h};v)\\
&+\big(\sigma^{n}-\sigma^{e}_{n},\psi\big)-
\big(\Delta_{(\beta-2)/2}(\varphi^{n}-\varphi^{e}_{n}),\psi\big)
+\big(\varphi^{n}-\varphi^{e}_{n},\phi\big)+\big(\psi^{n}-\psi^{e}_{n},\frac{\partial \phi}{\partial x}\big)\\
&-\sum_{k=1}^{K}\big(((\psi^{n}-\psi^{e}_{n})^{-}\phi^{-})_{k+\frac{1}{2}}
-((\psi^{n}-\psi^{e}_{n})^{-}\phi^{+})_{k-\frac{1}{2}}\big)
+\big(\psi^{n}-\psi^{e}_{n},\eta\big)+\big(\pi^{n}-\pi^{e}_{n},\frac{\partial \eta}{\partial x}\big)\\
&-\sum_{k=1}^{K}\big(((\pi^{n}-\pi^{e}_{n})^{+}\eta^{-})_{k+\frac{1}{2}}
-((\pi^{n}-\pi^{e}_{n})^{+}\eta^{+})_{k-\frac{1}{2}}\big)=0.\\
\end{split}
\end{equation}

Following the proof of Theorem \ref{th1}, we take the test functions
\begin{equation}\label{91h}
\begin{split}
v=\pi^{n},\quad  \psi=-\varphi^{n}+\sigma^{n},\quad\phi=\pi^{n},\quad\eta=\psi^{n},
\end{split}
\end{equation}
 we obtain
\begin{equation}\label{1dv1ggh}
\begin{split}
\bigg(\sum_{j=1}^{S}&W(\alpha_{j})\Delta\tau_{j}\delta_{t}^{\alpha_{j}}(\pi^{n}-\pi^{e}_{n})
,\pi^{n}\bigg)
-\varepsilon\big(\sigma^{n}-\sigma^{e}_{n},\pi^{n}\big)
+\big(\gamma(x)^{n},\pi^{n}\big)-\sum_{k=1}^{K}\mathcal{H}_{k}(f;u,u_{h};\pi^{n})\\
&+\big(\sigma^{n}-\sigma^{e}_{n},-\varphi^{n}+\sigma^{n}\big)-
\big(\Delta_{(\beta-2)/2}(\varphi^{n}-\varphi^{e}_{n}),-\varphi^{n}+\sigma^{n}\big)
+\big(\varphi^{n}-\varphi^{e}_{n},\pi^{n}\big)+\big(\psi^{n}-\psi^{e}_{n},\frac{\partial \pi^{n}}{\partial x}\big)\\
&-\sum_{k=1}^{K}\big(((\psi^{n}-\psi^{e}_{n})^{-}(\pi^{n})^{-})_{k+\frac{1}{2}}
-((\psi^{n}-\psi^{e}_{n})^{-}(\pi^{n})^{+})_{k-\frac{1}{2}}\big)
+\big(\psi^{n}-\psi^{e}_{n},\psi^{n}\big)+\big(\pi^{n}-\pi^{e}_{n},\frac{\partial \psi^{n}}{\partial x}\big)\\
&-\sum_{k=1}^{K}\big(((\pi^{n}-\pi^{e}_{n})^{+}(\psi^{n})^{-})_{k+\frac{1}{2}}
-((\pi^{n}-\pi^{e}_{n})^{+}(\psi^{n})^{+})_{k-\frac{1}{2}}\big)=0,\\
\end{split}
\end{equation}
by the properties of the projection $P^{+}$ and $P^{-}$, we obtain
\begin{equation}\label{1dv1ff}
\begin{split}
\bigg(\sum_{j=1}^{S}&W(\alpha_{j})\Delta\tau_{j}\delta_{t}^{\alpha_{j}}
(\pi^{n}-\pi^{e}_{n}),\pi^{n}\bigg)-\varepsilon\big(\sigma^{n}-\sigma^{e}_{n},\pi^{n}\big)
+\big(\gamma(x)^{n},\pi^{n}\big)-\sum_{k=1}^{K}\mathcal{H}_{k}(f;u,u_{h};\pi^{n})\\
&+\big(\sigma^{n}-\sigma^{e}_{n},-\varphi^{n}+\sigma^{n}\big)-
\big(\Delta_{(\beta-2)/2}(\varphi^{n}-\varphi^{e}_{n}),-\varphi^{n}+\sigma^{n}\big)
+\big(\varphi^{n}-\varphi^{e}_{n},\pi^{n}\big)\\
&+\sum_{k=1}^{K}\big(((\psi^{e}_{n})^{-}(\pi^{n})^{-})_{k+\frac{1}{2}}
-((\psi^{e}_{n})^{-}(\pi^{n})^{+})_{k-\frac{1}{2}}\big)
+\big(\psi^{n}-\psi^{e}_{n},\psi^{n}\big)\\
&+\sum_{k=1}^{K}\big(((\pi^{e}_{n})^{+}(\psi^{n})^{-})_{k+\frac{1}{2}}
-((\pi^{e}_{n})^{+}(\psi^{n})^{+})_{k-\frac{1}{2}}\big)=0.\\
\end{split}
\end{equation}
Employing Young's inequality and  Lemma \ref{lga2} and  the interpolation property \eqref{sth}, \eqref{sthzz} and \eqref{sthzzc}, we obtain
\begin{equation}\label{1dv1}
\begin{split}
&\bigg(\sum_{j=1}^{S}W(\alpha_{j})\Delta\tau_{j}\delta_{t}^{\alpha_{j}}\pi^{n},\pi^{n}\bigg)
+\big(\sigma^{n},\sigma^{n}\big)+
\big(\Delta_{(\beta-2)/2}\varphi^{n},\varphi^{n}\big)
+\big(\psi^{n},\psi^{n}\big)\\
&\leq
C\big(h^{2N+2}+(\Delta t)^{2+\theta}+\theta^{4}\big)+c_{2}\|\sigma^{n}\|^{2}_{L^{2}(\Omega)}+ c_{3}\|\varphi^{n}\|^{2}_{L^{2}(\Omega)}+c\|\pi^{n}\|^{2}_{L^{2}(\Omega)}
+c_{1}\|\psi^{n}\|^{2}_{L^{2}(\Omega)}\\
&\quad-\frac{1}{4}\kappa(f^{*};u_{h}^{n})
+(C+C_{*}(\|\pi^{n}\|_{\infty}+h^{-1}\|e_{u}\|_{\infty}^{2}))\|\pi^{n}\|^{2}
+(C+C_{*}h^{-1}\|e_{u}\|_{\infty}^{2})h^{2N+1}.
\end{split}
\end{equation}
Recalling Lemma \ref{lg} and provided $c_{i},\,\,i=1,2$ are sufficiently small such that $c_{i}\leq1$, we obtain
\begin{equation}\label{1dv1}
\begin{split}
&\bigg(\sum_{j=1}^{S}W(\alpha_{j})\Delta\tau_{j}\delta_{t}^{\alpha_{j}}\pi^{n},\pi^{n}\bigg)
\leq
C\big(h^{2N+1}+(\Delta t)^{2+\theta}+\theta^{4}\big)+c_{4}\|\pi^{n}\|^{2}_{L^{2}(\Omega)}.
\end{split}
\end{equation}
It then follows that
\begin{equation}\label{1dv1}
\begin{split}
\bigg(\sum_{j=1}^{S}
\frac{W(\alpha_{j})\Delta\tau_{j}}{\lambda_{j}}\pi^{n},\pi^{n}_{h}\bigg)\leq&
\bigg(\sum_{j=1}^{S}
\frac{W(\alpha_{j})\Delta\tau_{j}}{\lambda_{j}}\sum_{l=1}^{n-1}
(a_{n-l-1}^{\alpha_{j}}-a_{n-l}^{\alpha_{j}})\pi^{l},\pi^{n}\bigg)\\
&+
\bigg(\sum_{j=1}^{S}
\frac{W(\alpha_{j})\Delta\tau_{j}}{\lambda_{j}}a_{n-1}^{\alpha_{j}}\pi^{0},\pi^{n}\bigg)+ c\|\pi^{n}\|^{2}_{L^{2}(\Omega)}+C\big(h^{2N+1}+(\Delta t)^{2+\theta}+\theta^{4}\big).\\
\end{split}
\end{equation}
Employing Young's inequality, we obtain
\begin{equation}\label{1dv1}
\begin{split}
\|\pi^{n}\|^{2}_{L^{2}(\Omega_{h})}\leq &\sum_{j=1}^{S}
\frac{W(\alpha_{j})\Delta\tau_{j}}{\lambda_{j}}Q\sum_{l=1}^{n-1}
(a_{n-l-1}^{\alpha_{j}}-a_{n-l}^{\alpha_{j}})\|\pi^{l}\|^{2}_{L^{2}(\Omega_{h})}
+(cQ+\frac{1}{4})\sum_{j=1}^{S}
\frac{W(\alpha_{j})\Delta\tau_{j}}{\lambda_{j}}Q\|\pi^{n}\|^{2}_{L^{2}(\Omega_{h})}\\
&+C\sum_{j=1}^{S}
\frac{W(\alpha_{j})\Delta\tau_{j}}{\lambda_{j}}Qa_{n-1}^{\alpha_{j}}h^{2N+2}+C\sum_{j=1}^{S}
\frac{W(\alpha_{j})\Delta\tau_{j}}{\lambda_{j}}Qa_{n-1}^{\alpha_{j}}\big(h^{2N+1}+(\Delta t)^{2+\theta}+\theta^{4}\big),\\
\end{split}
\end{equation}
  provided $c$ is sufficiently small such that $\frac{3}{4}-cQ>0$, we obtain that
\begin{equation}\label{1dv1qw}
  \begin{split}
&\|\pi^{n}\|^{2}_{L^{2}(\Omega)}\leq C\bigg(\sum_{j=1}^{S}
\frac{W(\alpha_{j})\Delta\tau_{j}}{\lambda_{j}}Q\sum_{l=1}^{n-1}
(a_{n-l-1}^{\alpha_{j}}-a_{n-l}^{\alpha_{j}})\|\pi^{l}\|^{2}_{L^{2}(\Omega)}
+\sum_{j=1}^{S}
\frac{W(\alpha_{j})\Delta\tau_{j}}{\lambda_{j}}Qa_{n-1}^{\alpha_{j}}\big(h^{2N+1}+(\Delta t)^{2+\theta}+\theta^{4}\big)\bigg).\\
\end{split}
\end{equation}
Obviously the theorem holds for $n = 0$. Assume that it is valid for $n = 1, 2, ..., m-1$. Then,
by \eqref{1dv1qw}, we have
\begin{equation}\label{1dv1}
  \begin{split}
&\|\pi^{m}\|^{2}_{L^{2}(\Omega)}\leq C\bigg(\sum_{j=1}^{S}
\frac{W(\alpha_{j})\Delta\tau_{j}}{\lambda_{j}}Q\sum_{l=1}^{m-1}
(a_{n-l-1}^{\alpha_{j}}-a_{n-l}^{\alpha_{j}})\|\pi^{l}\|^{2}_{L^{2}(\Omega)}
+\sum_{j=1}^{S}
\frac{W(\alpha_{j})\Delta\tau_{j}}{\lambda_{j}}Qa_{n-1}^{\alpha_{j}}\big(h^{2N+1}+(\Delta t)^{2+\theta}+\theta^{4}\big)\bigg)\\
&\leq C\bigg(\sum_{j=1}^{S}
\frac{W(\alpha_{j})\Delta\tau_{j}}{\lambda_{j}}Q\sum_{l=1}^{m-1}
(a_{n-l-1}^{\alpha_{j}}-a_{n-l}^{\alpha_{j}})\big(h^{2N+1}+(\Delta t)^{2+\theta}+\theta^{4}\big)
+\sum_{j=1}^{S}
\frac{W(\alpha_{j})\Delta\tau_{j}}{\lambda_{j}}Qa_{n-1}^{\alpha_{j}}\big(h^{2N+1}+(\Delta t)^{4+\theta}+\theta^{4}\big)\bigg)
\\
&= C\big(h^{2N+1}+(\Delta t)^{2+\theta}+\theta^{4}\big).
\end{split}
\end{equation}
Finally, by using triangle inequality and standard approximation theory, we can get \eqref{tt7mk}. $\quad\Box$
\section{ LDG method for  the nonlinear distributed-order time and
space-fractional   Schr\"{o}dinger type equations}\label{s4}
\subsection{ LDG method for  the nonlinear distributed-order time and
space-fractional   Schr\"{o}dinger  equation}\label{s41}
 We rewrite the
fractional derivative as a composite of first order derivatives and a fractional integral
to recover the equation to a low order system. However, for the first order system,
alternating fluxes are used. We introduce three variables $e,r,s$ and set
\begin{equation}\label{1a}
\begin{split}
&e=\Delta_{(\beta-2)/2}r, \quad r=\frac{\partial}{\partial x}s,\quad s=\frac{\partial}{\partial x}u,
\end{split}
\end{equation}
then, the nonlinear distributed-order time and
space-fractional  Schr\"{o}dinger  problem can be rewritten as
\begin{equation}\label{1bchcf}
\begin{split}
&i\mathcal{D}_{t}^{W(\alpha)}{}u+\varepsilon_{1}e+ \varepsilon_{2}f(|u|^{2})u=0,\\
&e=\Delta_{(\beta-2)/2}r,\quad r=\frac{\partial}{\partial x}s,\quad s=\frac{\partial}{\partial x}u.\\
\end{split}
\end{equation}
For actual numerical implementation, it might be more efficient if we decompose the complex function
$u(x,t)$ into its real and imaginary parts by writing
\begin{equation}\label{1b}
\begin{split}
u(x,t)=p(x,t)+iq(x,t),
\end{split}
\end{equation}
where $p$, $q$ are real functions. Under the new notation, the problem \eqref{1bchcf} can be written as
\begin{equation}\label{1bx}
\begin{split}
&\mathcal{D}_{t}^{W(\alpha)}{}p+\varepsilon_{1}e+ \varepsilon_{2}f(p^{2}+q^{2})q=0,\\
&e=\Delta_{(\beta-2)/2}r,\quad r=\frac{\partial}{\partial x}s,\quad s=\frac{\partial}{\partial x}q,\\
&\mathcal{D}_{t}^{W(\alpha)}{}q-\varepsilon_{1}l- \varepsilon_{2}f(p^{2}+q^{2})p=0,\\
&l=\Delta_{(\alpha-2)/2}w,\quad w=\frac{\partial}{\partial x}z,\quad z=\frac{\partial}{\partial x}p.\\
\end{split}
\end{equation}
 Let $p_{h}^{n}, q_{h}^{n}, e_{h}^{n},l_{h}^{n},r_{h}^{n},s_{h}^{n},w_{h}^{n},z_{h}^{n}\in V_{k}^{N}$ be the approximation of
$p(.,t_{n}),q(.,t_{n}),e(.,t_{n}),l(.,t_{n})r(.,t_{n}),s(.,t_{n}),$ \\
$w(.,t_{n}),z(.,t_{n})$ respectively.
 We a fully discrete local discontinuous Galerkin scheme as follows: find $p_{h}^{n}, q_{h}^{n}, e_{h}^{n},$\\
 $l_{h}^{n},r_{h}^{n},s_{h}^{n},w_{h}^{n},z_{h}^{n}\in V_{k}^{N}$, such that for all test functions $\vartheta_{1},\rho,\phi,\varphi,\chi,\varrho,\psi,\zeta\in V_{k}^{N}$,
\begin{equation}\label{1bch1nb}
\begin{split}
&\bigg(\sum_{j=1}^{S}W(\alpha_{j})\Delta\tau_{j}\delta_{t}^{\alpha_{j}}p_{h}^{n},\vartheta_{1}\bigg)_{D^{k}}
+\varepsilon_{1}\big(e_{h}^{n},\vartheta_{1}\big)_{D^{k}}+  \varepsilon_{2}\big(f((p_{h}^{n})^{2}+(q_{h}^{n})^{2})q_{h}^{n},\vartheta_{1}\big)_{D^{k}}=0,\\
&\big(e_{h}^{n},\rho\big)_{D^{k}}=\big(\Delta_{(\beta-2)/2}r_{h}^{n},\rho\big)_{D^{k}},\\
&\big(r_{h}^{n},\phi\big)_{D^{k}}=\big(\frac{\partial}{\partial x}s_{h}^{n},\phi\big)_{D^{k}},\\
&\big(s_{h}^{n},\varphi\big)_{D^{k}}=\big(\frac{\partial}{\partial x}q_{h}^{n},\varphi\big)_{D^{k}},\\
&\big(\sum_{j=1}^{S}W(\alpha_{j})\Delta\tau_{j}\delta_{t}^{\alpha_{j}}q_{h}^{n},\chi\big)_{D^{k}}
-\varepsilon_{1}\big(l_{h}^{n},\chi\big)_{D^{k}}- \varepsilon_{2}\big(f((p_{h}^{n})^{2}+(q_{h}^{n})^{2}))p_{h}^{n},\chi\big)_{D^{k}}=0,\\
&\big(l_{h}^{n},\varrho\big)_{D^{k}}=\big(\Delta_{(\beta-2)/2}w_{h}^{n},\varrho\big)_{D^{k}},\\
&\big(w_{h}^{n},\psi\big)_{D^{k}}=\big(\frac{\partial}{\partial x}z_{h}^{n},\psi\big)_{D^{k}},\\
&\big(z_{h}^{n},\zeta\big)_{D^{k}}=\big(\frac{\partial}{\partial x}p_{h}^{n},\zeta\big)_{D^{k}}.\\
\end{split}
\end{equation}
Applying integration by parts to \eqref{1bch1nb}, and replacing the fluxes at the interfaces
by the corresponding numerical fluxes, we obtain
\begin{equation}\label{1bch2}
\begin{split}
&\big(\sum_{j=1}^{S}W(\alpha_{j})\Delta\tau_{j}\delta_{t}^{\alpha_{j}}p_{h}^{n},\vartheta_{1}\big)_{D^{k}}
+\varepsilon_{1}\big(e_{h}^{n},\vartheta_{1}\big)_{D^{k}}+ \varepsilon_{2}\big(\big(f((p_{h}^{n})^{2}+(q_{h}^{n})^{2}))q_{h}^{n},\vartheta_{1}\big)_{D^{k}}=0,\\
&\big(e_{h}^{n},\rho\big)_{D^{k}}=\big(\Delta_{(\beta-2)/2}r_{h}^{n},\rho\big)_{D^{k}},\\
&\big(r_{h}^{n},\phi\big)_{D^{k}}=-\big(s_{h}^{n},\phi_{x}\big)_{D^{k}}
+\big((\widehat{s}^{n}_{h}\phi^{-})_{k+\frac{1}{2}}
-(\widehat{s}^{n}_{h}\phi^{+})_{k-\frac{1}{2}}\big),\\
&\big(s_{h}^{n},\varphi\big)_{D^{k}}=-\big(q_{h}^{n},\varphi_{x}\big)_{D^{k}}
+\big((\widehat{q}^{n}_{h}\varphi^{-})_{k+\frac{1}{2}}
-(\widehat{q}^{n}_{h}\varphi^{+})_{k-\frac{1}{2}}\big),\\
&\big(\sum_{j=1}^{S}W(\alpha_{j})\Delta\tau_{j}\delta_{t}^{\alpha_{j}}q_{h}^{n},\chi\big)_{D^{k}}
-\varepsilon_{1}\big(l_{h}^{n},\chi\big)_{D^{k}}- \varepsilon_{2}\big(\big(f((p_{h}^{n})^{2}+(q_{h}^{n})^{2}))p_{h}^{n},\chi\big)_{D^{k}}=0,\\
&\big(l_{h}^{n},\varrho\big)_{D^{k}}=\big(\Delta_{(\beta-2)/2}w_{h}^{n},\varrho\big)_{D^{k}},\\
&\big(w_{h}^{n},\psi\big)_{D^{k}}=-\big(z_{h}^{n},\psi_{x}\big)_{D^{k}}+
\big((\widehat{z}^{n}_{h}\psi^{-})_{k+\frac{1}{2}}
-(\widehat{z}^{n}_{h}\psi^{+})_{k-\frac{1}{2}}\big),\\
&\big(z_{h}^{n},\zeta\big)_{D^{k}}=-\big(p_{h}^{n},\zeta_{x}\big)_{D^{k}}+\big
((\widehat{p}^{n}_{h}\zeta^{-})_{k+\frac{1}{2}}
-(\widehat{p}^{n}_{h}\zeta^{+})_{k-\frac{1}{2}}\big).\\
\end{split}
\end{equation}
The numerical traces $( \widehat{p}^{n}_{h}, \widehat{q}^{n}_{h},\widehat{s}^{n}_{h},\widehat{z}^{n}_{h})$ are defined on interelement faces as the
alternating fluxes \cite{doi:10.1137/S0036142997316712,ref1}:
\begin{equation}\label{fl}
\begin{split}
\widehat{p}^{n}_{h}=(p^{n}_{h})^{-},\,\,\widehat{s}^{n}_{h}
=(s^{n}_{h})^{+},\,\,\,\,
\widehat{q}^{n}_{h}=(q^{n}_{h})^{-},\,\,
\widehat{z}^{n}_{h}=(z^{n}_{h})^{+}.
\end{split}
\end{equation}
Note that we can also choose
\begin{equation}\label{1a}
\begin{split}
\widehat{p}^{n}_{h}=(p^{n}_{h})^{+},\,\,\widehat{s}^{n}_{h}
=(s^{n}_{h})^{-},\,\,\,\,
\widehat{q}^{n}_{h}=(q^{n}_{h})^{+},\,\,
\widehat{z}^{n}_{h}=(z^{n}_{h})^{-}.
\end{split}
\end{equation}
\subsubsection{ The analysis of stability for fully discrete scheme }\label{s42}  In order to carry out the analysis of the LDG scheme, we have the following results.
\begin{thm}\label{tt4h} ($L^{2}$ stability).
The semidiscrete scheme \eqref{1bch2} is   stable, and
\begin{equation}\label{1bch1nn}
\begin{split}
&\|p_{h}^{n}\|^{2}_{L^{2}(\Omega)}+ \|q_{h}^{n}\|^{2}_{L^{2}(\Omega)}\leq
C\big(\|p_{h}^{0}\|^{2}_{L^{2}(\Omega)}+\|q_{h}^{0}\|^{2}_{L^{2}(\Omega)}\big). 
\end{split}
\end{equation}

\end{thm}
\textbf{Proof.}
Set $(\vartheta_{1},\rho,\phi,\varphi,\chi,\varrho,\psi,\zeta)=(p_{h}^{n},-r_{h}^{n}+e_{h}^{n},p_{h}^{n}
,-z_{h}^{n},q_{h}^{n},l_{h}^{n}-w_{h}^{n},-q_{h}^{n},s_{h}^{n})$ in \eqref{1bch2}, we get
\begin{equation}\label{1bch1bb}
\begin{split}
&\big(\sum_{j=1}^{S}W(\alpha_{j})\Delta\tau_{j}\delta_{t}^{\alpha_{j}}p_{h}^{n},p_{h}^{n}\big)_{D^{k}}
+\big(\sum_{j=1}^{S}W(\alpha_{j})\Delta\tau_{j}\delta_{t}^{\alpha_{j}}q_{h}^{n},q_{h}^{n}\big)_{D^{k}}
+\big(e_{h}^{n},e_{h}^{n}\big)_{D^{k}}+\big(l_{h}^{n},l_{h}^{n}\big)_{D^{k}}
+\big(\Delta_{(\beta-2)/2}w_{h}^{n},w_{h}^{n}\big)_{D^{k}}\\
&\quad+\big(\Delta_{(\beta-2)/2}r_{h}^{n},r_{h}^{n}\big)_{D^{k}}=\big(\Delta_{(\beta-2)/2}w_{h}^{n},l_{h}^{n}\big)_{D^{k}}
+\big(\Delta_{(\beta-2)/2}r_{h}^{n},e_{h}^{n}\big)_{D^{k}}-\big( r_{h}^{n},p_{h}^{n}\big)_{D^{k}}+\big( w_{h}^{n},q_{h}^{n}\big)_{D^{k}}+\big( e_{h}^{n},r_{h}^{n}\big)_{D^{k}}\\
&\quad-\big(s_{h}^{n},(p_{h}^{n})_{x}\big)_{D^{k}}+\big(q_{h}^{n},(z_{h}^{n})_{x}\big)_{D^{k}}
+\big(z_{h}^{n},(q_{h}^{n})_{x}\big)_{D^{k}}
-\big(p_{h}^{n},(s_{h}^{n})_{x}\big)_{D^{k}}+\big( l_{h}^{n},w_{h}^{n}\big)_{D^{k}}-\varepsilon_{1}\big( e_{h}^{n},p_{h}^{n}\big)_{D^{k}}\\
&\quad+\varepsilon_{1}\big( l_{h}^{n},q_{h}^{n}\big)_{D^{k}}+\big((\widehat{s}^{n}_{h}(p_{h}^{n})^{-})_{k+\frac{1}{2}}
-(\widehat{s}^{n}_{h}(p_{h}^{n})^{+})_{k-\frac{1}{2}}\big)-
\big((\widehat{q}^{n}_{h}(z_{h}^{n})^{-})_{k+\frac{1}{2}}
-(\widehat{q}^{n}_{h}(z_{h}^{n})^{+})_{k-\frac{1}{2}}\big)\\
&\quad-\big((\widehat{z}^{n}_{h}(q_{h}^{n})^{-})_{k+\frac{1}{2}}
-(\widehat{z}^{n}_{h}(q_{h}^{n})^{+})_{k-\frac{1}{2}}\big)
+\big((\widehat{p}^{n}_{h}(s_{h}^{n})^{-})_{k+\frac{1}{2}}
-(\widehat{p}^{n}_{h}(s_{h}^{n})^{+})_{k-\frac{1}{2}}\big).\\
\end{split}
\end{equation}
Summing over $k$, with the definition \eqref{fl} of the numerical fluxes
and with simple algebraic manipulations, we easily obtain
\begin{equation}\label{1bch1}
\begin{split}
&\big(\sum_{j=1}^{S}W(\alpha_{j})\Delta\tau_{j}\delta_{t}^{\alpha_{j}}p_{h}^{n},p_{h}^{n}\big)
+\big(\sum_{j=1}^{S}W(\alpha_{j})\Delta\tau_{j}\delta_{t}^{\alpha_{j}}q_{h}^{n},q_{h}^{n}\big)
+\big(e_{h}^{n},e_{h}^{n}\big)+\big(l_{h}^{n},l_{h}^{n}\big)
+\big(\Delta_{(\beta-2)/2}w_{h}^{n},w_{h}^{n}\big)\\
&\quad+\big(\Delta_{(\beta-2)/2}r_{h}^{n},r_{h}^{n}\big)
=\big(\Delta_{(\beta-2)/2}w_{h}^{n},l_{h}^{n}\big)
+\big(\Delta_{(\beta-2)/2}r_{h}^{n},e_{h}^{n}\big)-\big( r_{h}^{n},p_{h}^{n}\big)+\big( w_{h}^{n},q_{h}^{n}\big)+\big( e_{h}^{n},r_{h}^{n}\big).
\end{split}
\end{equation}
Employing Young's inequality and  Lemma \ref{lga2}, we obtain
\begin{equation}\label{1bch1}
\begin{split}
&\bigg(\sum_{j=1}^{S}W(\alpha_{j})\Delta\tau_{j}\delta_{t}^{\alpha_{j}}p_{h}^{n},p_{h}^{n}\bigg)
+\bigg(\sum_{j=1}^{S}W(\alpha_{j})\Delta\tau_{j}\delta_{t}^{\alpha_{j}}q_{h}^{n},q_{h}^{n}\bigg)
+\big(e_{h}^{n},e_{h}^{n}\big)_{D^{k}}+\big(l_{h}^{n},l_{h}^{n}\big)
+\big(\Delta_{(\beta-2)/2}w_{h}^{n},w_{h}^{n}\big)\\
&\quad+\big(\Delta_{(\beta-2)/2}r_{h}^{n},r_{h}^{n}\big)_{\Omega_{h}}\leq c\|p_{h}^{n}\|^{2}_{L^{2}(\Omega)}+c\|q_{h}^{n}\|^{2}_{L^{2}(\Omega)}
+c_{6}\|w_{h}^{n}\|^{2}_{L^{2}(\Omega)}+c_{5}\|r_{h}^{n}\|^{2}_{L^{2}(\Omega)}
+c_{1}\|e_{h}^{n}\|^{2}_{L^{2}(\Omega)}
\\
&\qquad\qquad\qquad\qquad\qquad\qquad\qquad+c_{2}\|l_{h}^{n}\|^{2}_{L^{2}(\Omega)}.\\
\end{split}
\end{equation}
Recalling Lemma \ref{lg} and provided $c_{i},\,\,i=1,2$ are sufficiently small such that $c_{i}\leq1$, we obtain that
\begin{equation}\label{1bch1n}
\begin{split}
&\bigg(\sum_{j=1}^{S}W(\alpha_{j})\Delta\tau_{j}\delta_{t}^{\alpha_{j}}p_{h}^{n},p_{h}^{n}\bigg)
+\bigg(\sum_{j=1}^{S}W(\alpha_{j})\Delta\tau_{j}\delta_{t}^{\alpha_{j}}q_{h}^{n},q_{h}^{n}\bigg)
\leq c(\|p_{h}^{n}\|^{2}_{L^{2}(\Omega)}+\|q_{h}^{n}\|^{2}_{L^{2}(\Omega)}).\\
\end{split}
\end{equation}
It then follows that
\begin{equation}\label{1bch1}
\begin{split}
&\bigg(\sum_{j=1}^{S}
\frac{W(\alpha_{j})\Delta\tau_{j}}{\lambda_{j}}p^{n},p^{n}_{h}\bigg)
+
\bigg(\sum_{j=1}^{S}
\frac{W(\alpha_{j})\Delta\tau_{j}}{\lambda_{j}}q^{n},q^{n}_{h}\bigg)\\
&\leq \bigg(\sum_{j=1}^{S}
\frac{W(\alpha_{j})\Delta\tau_{j}}{\lambda_{j}}\sum_{l=1}^{n-1}
(a_{n-l-1}^{\alpha_{j}}-a_{n-l}^{\alpha_{j}})p^{l}_{h},p^{n}_{h}\bigg)+
\bigg(\sum_{j=1}^{S}
\frac{W(\alpha_{j})\Delta\tau_{j}}{\lambda_{j}}a_{n-1}^{\alpha_{j}}p^{0}_{h},p^{n}_{h}\bigg)\\
&\quad+
\bigg(\sum_{j=1}^{S}
\frac{W(\alpha_{j})\Delta\tau_{j}}{\lambda_{j}}\sum_{l=1}^{n-1}
(a_{n-l-1}^{\alpha_{j}}-a_{n-l}^{\alpha_{j}})q^{l}_{h},q^{n}_{h}\bigg)
+
\bigg(\sum_{j=1}^{S}
\frac{W(\alpha_{j})\Delta\tau_{j}}{\lambda_{j}}a_{n-1}^{\alpha_{j}}q^{0}_{h},q^{n}_{h}\bigg) +c(\|p_{h}^{n}\|^{2}_{L^{2}(\Omega)}+\|q_{h}^{n}\|^{2}_{L^{2}(\Omega)}).
\end{split}
\end{equation}
Employing Young's inequality, we obtain
\begin{equation}\label{1bch1}
\begin{split}
&\|p_{h}^{n}\|^{2}_{L^{2}(\Omega)}
+
\|q_{h}^{n}\|^{2}_{L^{2}(\Omega)}\\
&\leq \frac{1}{2}\sum_{j=1}^{S}
\frac{W(\alpha_{j})\Delta\tau_{j}}{\lambda_{j}}Q\bigg(\sum_{l=1}^{n-1}
(a_{n-l-1}^{\alpha_{j}}-a_{n-l}^{\alpha_{j}})\big(\|p_{h}^{l}\|^{2}_{L^{2}(\Omega)}
+\|p_{h}^{n}\|^{2}_{L^{2}(\Omega)}\big)\bigg)+
\frac{1}{2}\sum_{j=1}^{S}
\frac{W(\alpha_{j})\Delta\tau_{j}}{\lambda_{j}}Qa_{n-1}^{\alpha_{j}}
(\|p_{h}^{0}\|^{2}_{L^{2}(\Omega)}\\
&\quad+\|p_{h}^{n}\|^{2}_{L^{2}(\Omega)})+
\frac{1}{2}\sum_{j=1}^{S}
\frac{W(\alpha_{j})\Delta\tau_{j}}{\lambda_{j}}Q\bigg(\sum_{l=1}^{n-1}
(a_{n-l-1}^{\alpha_{j}}-a_{n-l}^{\alpha_{j}})\big(\|q_{h}^{l}\|^{2}_{L^{2}(\Omega)}
+\|q_{h}^{n}\|^{2}_{L^{2}(\Omega)}\big)\bigg)\\
&\quad+\frac{1}{2}
\sum_{j=1}^{S}
\frac{W(\alpha_{j})\Delta\tau_{j}}{\lambda_{j}}Qa_{n-1}^{\alpha_{j}}\big(\|q_{h}^{0}\|^{2}_{L^{2}(\Omega)}
+\|q_{h}^{n}\|^{2}_{L^{2}(\Omega)}\big) +cQ(\|p_{h}^{n}\|^{2}_{L^{2}(\Omega)}+\|q_{h}^{n}\|^{2}_{L^{2}(\Omega)}),
\end{split}
\end{equation}
provided $c$ is sufficiently small such that $\frac{1}{2}-cQ>0$, we obtain that
\begin{equation}\label{1bch1jk}
\begin{split}
&\|p_{h}^{n}\|^{2}_{L^{2}(\Omega)}
+
\|q_{h}^{n}\|^{2}_{L^{2}(\Omega)}\\
&\leq
 C\sum_{j=1}^{S}
\frac{W(\alpha_{j})\Delta\tau_{j}}{\lambda_{j}}Q\bigg(\sum_{l=1}^{n-1}
(a_{n-l-1}^{\alpha_{j}}-a_{n-l}^{\alpha_{j}})\big(\|p_{h}^{l}\|^{2}_{L^{2}(\Omega)}
+\|q_{h}^{l}\|^{2}_{L^{2}(\Omega)}\big)\bigg)+
C\sum_{j=1}^{S}
\frac{W(\alpha_{j})\Delta\tau_{j}}{\lambda_{j}}Qa_{n-1}^{\alpha_{j}}
\big(\|p_{h}^{0}\|^{2}_{L^{2}(\Omega)}\\
&\quad+\|q_{h}^{0}\|^{2}_{L^{2}(\Omega)}\big) . 
\end{split}
\end{equation}
Obviously the theorem holds for $n = 0$. Assume that it is valid for $n = 1, 2, ..., m-1$. Then,
by \eqref{1bch1jk}, we have
\begin{equation}\label{1bch1}
\begin{split}
&\|p_{h}^{m}\|^{2}_{L^{2}(\Omega)}
+
\|q_{h}^{m}\|^{2}_{L^{2}(\Omega)}\\
&\leq C\sum_{j=1}^{S}
\frac{W(\alpha_{j})\Delta\tau_{j}}{\lambda_{j}}Q\bigg(\sum_{l=1}^{m-1}
(a_{n-l-1}^{\alpha_{j}}-a_{n-l}^{\alpha_{j}})\big(\|p_{h}^{l}\|^{2}_{L^{2}(\Omega)}
+\|q_{h}^{l}\|^{2}_{L^{2}(\Omega)}\big)\bigg)+
C\sum_{j=1}^{S}
\frac{W(\alpha_{j})\Delta\tau_{j}}{\lambda_{j}}Qa_{n-1}^{\alpha_{j}}
\big(\|p_{h}^{0}\|^{2}_{L^{2}(\Omega)}\\
&\quad\quad\quad\qquad+\|q_{h}^{0}\|^{2}_{L^{2}(\Omega)}\big).\\
&\leq C\sum_{j=1}^{S}
\frac{W(\alpha_{j})\Delta\tau_{j}}{\lambda_{j}}Q\bigg(\sum_{l=1}^{m-1}
(a_{n-l-1}^{\alpha_{j}}-a_{n-l}^{\alpha_{j}})\big(\|p_{h}^{0}\|^{2}_{L^{2}(\Omega)}
+\|q_{h}^{0}\|^{2}_{L^{2}(\Omega)}
\big)\bigg)+
C\sum_{j=1}^{S}
\frac{W(\alpha_{j})\Delta\tau_{j}}{\lambda_{j}}Qa_{n-1}^{\alpha_{j}}
\big(\|p_{h}^{0}\|^{2}_{L^{2}(\Omega)}\\
&\quad\quad\quad\qquad+\|q_{h}^{0}\|^{2}_{L^{2}(\Omega)}\big).\\
&=C\big(\|p_{h}^{0}\|^{2}_{L^{2}(\Omega)}+\|q_{h}^{0}\|^{2}_{L^{2}(\Omega)}\big). \quad\Box\\
\end{split}
\end{equation}

\subsubsection{Error estimates} \label{s43}
We consider the linear distributed-order time and space-fractional Schr\"{o}dinger equation
\begin{equation}\label{1bchj}
\begin{split}
&i\mathcal{D}_{t}^{W(\alpha)}{}u-\varepsilon_{1}(-\Delta) ^{\frac{\alpha}{2}}u+\varepsilon_{2}u=0.\\
\end{split}
\end{equation}
It is easy to verify that the exact solution of the above  \eqref{1bchj} satisfies
\begin{equation}\label{1bch2j}
\begin{split}
&\bigg( \sum_{j=1}^{S}W(\alpha_{j})\Delta\tau_{j}\delta_{t}^{\alpha_{j}}p^{n},\vartheta_{1}\bigg)_{D^{k}}
+\varepsilon_{1}\big(e^{n},\vartheta_{1}\big)_{D^{k}}+ \varepsilon_{2}\big(q^{n},\vartheta_{1}\big)_{D^{k}}+\big(\gamma(x)^{n},\vartheta_{1}\big)_{D^{k}}=0,\\
&\big(e^{n},\rho\big)_{D^{k}}=\big(\Delta_{(\beta-2)/2}r^{n},\rho\big)_{D^{k}},\\
&\big(r^{n},\phi\big)_{D^{k}}=-\big(s^{n},\phi_{x}\big)_{D^{k}}+\big((
\widehat{s}^{n}\phi^{-})_{k+\frac{1}{2}}
-(\widehat{s}^{n}\phi^{+})_{k-\frac{1}{2}}\big),\\
&\big(s^{n},\varphi\big)_{D^{k}}=-\big(q^{n},\varphi_{x}\big)_{D^{k}}
+\big((\widehat{q}^{n}\varphi^{-})_{k+\frac{1}{2}}
-(\widehat{q}^{n}\varphi^{+})_{k-\frac{1}{2}}\big),\\
&\bigg(\sum_{j=1}^{S}W(\alpha_{j})\Delta\tau_{j}\delta_{t}^{\alpha_{j}}q^{n},\chi\bigg)_{D^{k}}
-\varepsilon_{1}\big(l^{n},\chi\big)_{D^{k}}-\varepsilon_{2} \big(p^{n},\chi\big)_{D^{k}}=0,\\
&\big(l^{n},\varrho\big)_{D^{k}}=\big(\Delta_{(\beta-2)/2}w^{n},\varrho\big)_{D^{k}},\\
&\big(w^{n},\psi\big)_{D^{k}}=-\big(z^{n},\psi_{x}\big)_{D^{k}}+\big((
\widehat{z}^{n}\psi^{-})_{k+\frac{1}{2}}
-(\widehat{z}^{n}\psi^{+})_{k-\frac{1}{2}}\big),\\
&\big(z^{n},\zeta\big)_{D^{k}}=-\big(p^{n},\zeta_{x}\big)_{D^{k}}+\big((\widehat{p}^{n}\psi^{-})_{k+\frac{1}{2}}
-(\widehat{p}^{n}\zeta^{+})_{k-\frac{1}{2}}\big).\\
\end{split}
\end{equation}
Subtracting  \eqref{1bch2}  from  \eqref{1bch2j}, we can obtain the error equation
\begin{equation}\label{1bch2k2}
\begin{split}
\bigg( &\sum_{j=1}^{S}W(\alpha_{j})\Delta\tau_{j}\delta_{t}^{\alpha_{j}}(p^{n}-p_{h}^{n}),\vartheta_{1}\bigg)_{D^{k}}
+\big(\sum_{j=1}^{S}W(\alpha_{j})\Delta\tau_{j}\delta_{t}^{\alpha_{j}}(q^{n}-q_{h}^{n}),\chi\big)_{D^{k}}
-\big(\Delta_{(\beta-2)/2}(r^{n}-r_{h}^{n}),\rho\big)_{D^{k}}
\\
&-\big(\Delta_{(\beta-2)/2}(w^{n}-w_{h}^{n}),\varrho\big)_{D^{k}}+\big(s^{n}-s_{h}^{n},\phi_{x}\big)_{D^{k}} +\big(q^{n}-q_{h}^{n},\varphi_{x}\big)_{D^{k}}+\big(z^{n}-z_{h}^{n},\psi_{x}\big)_{D^{k}}
+\big(p^{n}-p_{h}^{n},\zeta_{x}\big)_{D^{k}}\\
&+\big(\gamma(x)^{n},\vartheta_{1}\big)_{D^{k}}+ \varepsilon_{2}\big(q^{n}-q_{h}^{n},\vartheta_{1}\big)_{D^{k}}- \varepsilon_{2}\big(p^{n}-p_{h}^{n},\chi\big)_{D^{k}}
+\big(r^{n}-r_{h}^{n},\phi\big)_{D^{k}}+\big(s^{n}-s_{h}^{n},\varphi\big)_{D^{k}}
+\big(l^{n}-l_{h}^{n},\varrho\big)_{D^{k}}\\
&+\big(e^{n}-e_{h}^{n},\rho\big)_{D^{k}}+\big(w^{n}-w_{h}^{n},\psi\big)_{D^{k}} +\big(z^{n}-z_{h}^{n},\zeta\big)_{D^{k}}-\varepsilon_{1}\big(l^{n}-l_{h}^{n},\chi\big)_{D^{k}}
+\varepsilon_{1}\big(e^{n}-e_{h}^{n},\vartheta_{1}\big)_{D^{k}}\\
&-\big(((\widehat{s}^{n}_{h}-\widehat{s}^{n})\phi^{-})_{k+\frac{1}{2}}
-((\widehat{s}^{n}_{h}-\widehat{s}^{n})\phi^{+})_{k-\frac{1}{2}}\big)
-\big(((\widehat{q}^{n}_{h}-\widehat{q}^{n})\varphi^{-})_{k+\frac{1}{2}}
-((\widehat{q}^{n}_{h}-\widehat{q}^{n})\varphi^{+})_{k-\frac{1}{2}}\big)\\
&-\big(((\widehat{z}^{n}_{h}-\widehat{z}^{n})\psi^{-})_{k+\frac{1}{2}}
-((\widehat{z}^{n}_{h}-\widehat{z}^{n})\psi^{+})_{k-\frac{1}{2}}\big)
-\big(((\widehat{p}^{n}_{h}-\widehat{p}^{n})\zeta^{-})_{k+\frac{1}{2}}
-((\widehat{p}^{n}_{h}-\widehat{p}^{n})\zeta^{+})_{k-\frac{1}{2}}\big)=0.\\
\end{split}
\end{equation}

Denoting
\begin{equation}\label{91h}
\begin{split}
&\pi^{n}=\mathcal{P}^{-}p^{n}-p_{h}^{n},\quad \pi^{e}_{n}=\mathcal{P}^{-}p^{n}-p^{n},\quad \epsilon^{n}=\mathcal{P}r^{n}-r_{h}^{n},\quad \epsilon^{e}_{n}=\mathcal{P}r^{n}-r^{n},\quad\phi^{n}=\mathcal{P}e^{n}-e_{h}^{n},\quad
\phi^{e}_{n}=\mathcal{P}e^{n}-e^{n},\\
& \tau^{n}=\mathcal{P}^{+}s^{n}-s_{h}^{n},\quad\tau^{e}_{n}=\mathcal{P}^{+}s^{n}-s^{n},\quad \sigma^{n}=\mathcal{P}^{-}q^{n}-q_{h}^{n}, \quad \sigma^{e}_{n}=\mathcal{P}^{-}q^{n}-q^{n},\quad\varpi^{n}=\mathcal{P}l^{n}-l_{h}^{n},\quad
\varpi^{e}_{n}=\mathcal{P}l^{n}-l^{n},  \\
&\varphi^{n}=\mathcal{P}w^{n}-w_{h}^{n},\quad \varphi^{e}_{n}=\mathcal{P}w^{n}-w^{n},\quad \vartheta^{n}=\mathcal{P}^{+}z^{n}-z_{h}^{n},\quad \vartheta^{e}_{n}=\mathcal{P}^{+}z^{n}-z^{n}.
\end{split}
\end{equation}

 \begin{lem}\label{lm1h}
\begin{equation}\label{91chvv}
\begin{split}
&\bigg(\sum_{j=1}^{S}W(\alpha_{j})\Delta\tau_{j}\delta_{t}^{\alpha_{j}}\pi^{n},\pi^{n}\bigg)
+\bigg(\sum_{j=1}^{S}W(\alpha_{j})
\Delta\tau_{j}\delta_{t}^{\alpha_{j}}\sigma^{n},\sigma^{n}\bigg)
+\big(\Delta_{(\beta-2)/2}\epsilon^{n},\epsilon^{n}\big)
+\big(\Delta_{(\beta-2)/2}\varphi^{n},\varphi^{n}\big)\\
&\quad+\big(\phi^{n},\phi^{n}\big)
+\big(\varpi^{n},\varpi^{n}\big)=Q_{1}+Q_{2}+Q_{3}+Q_{4}, \\
\end{split}
\end{equation}
where
\begin{subequations}\label{91chjm}
\begin{align}
&Q_{1}=-\big(\epsilon^{n},\pi^{n}\big)
+\big(\varphi^{n},\sigma^{n}\big)+\big(\Delta_{(\beta-2)/2}\epsilon^{n},\phi^{n}\big)
+\big(\Delta_{(\beta-2)/2}\varphi^{n},\varpi^{n}\big)\\
&\quad\quad\quad-\varepsilon_{1}\big(\phi^{n},\pi^{n}\big)
+\varepsilon_{1}\big(\varpi^{n},\sigma^{n}\big) +\big(\varpi^{n},\varphi^{n}\big)
+\big(\phi^{n},\epsilon^{n}\big),\\
&Q_{2}=\big(\tau^{e}_{n},\pi_{x}^{n}\big) -\big(\sigma^{e}_{n},\vartheta_{x}^{n}\big)-\big(\vartheta^{e}_{n},\sigma_{x}^{n}\big)
+\big(\pi^{e}_{n},\tau_{x}^{n}\big)
+\big(\vartheta^{e}_{n},\tau^{n}\big)-\big(\tau^{e}_{n},\vartheta^{n}\big),\\
&Q_{3}=\big( \sum_{j=1}^{S}W(\alpha_{j})\Delta\tau_{j}\delta_{t}^{\alpha_{j}}\pi^{e}_{n},\pi^{n}\big)
+\big(\sum_{j=1}^{S}W(\alpha_{j})\Delta\tau_{j}\delta_{t}^{\alpha_{j}}\sigma^{e}_{n},\sigma^{n}\big)
+\big(\varpi^{e},\varpi-\varphi\big)\\
&\quad\quad+\big(\phi_{n}^{e},\phi^{n}-\epsilon^{n}\big)
+\varepsilon_{2}\big(\sigma^{e}_{n},\pi^{n}\big)-\varepsilon_{2}\big(\pi^{e}_{n},\sigma^{n}\big)
+\big(\epsilon^{e}_{n},\pi^{n}\big)
-\big(\varphi^{e}_{n},\sigma^{n}\big)\\
&\quad\quad-\big(\Delta_{(\beta-2)/2}\epsilon^{e}_{n},\phi^{n}
-\epsilon^{n}\big)
-\big(\Delta_{(\beta-2)/2}\varphi^{e}_{n},\varpi^{n}-\varphi^{n}\big)
+\varepsilon_{1}\big(\phi_{n}^{e},\pi^{n}\big)
-\varepsilon_{1}\big(\varpi_{n}^{e},\sigma^{n}\big)\\
&\quad\quad+\big(\gamma(x)^{n},\pi^{n}\big)
,\\
&Q_{4}=-\sum_{k=1}^{K}\big(((\tau^{e}_{n})^{+}(\pi^{n})^{-})_{k+\frac{1}{2}}
-((\tau^{e}_{n})^{+}(\pi^{n})^{+})_{k-\frac{1}{2}}\big)
+\sum_{k=1}^{K}\big(((\sigma^{e}_{n})^{-}(\vartheta^{n})^{-})_{k+\frac{1}{2}}
-((\sigma^{e}_{n})^{-}(\vartheta^{n})^{+})_{k-\frac{1}{2}}\big)
\\
&\quad\quad\quad+\sum_{k=1}^{K}\big(((\vartheta^{e}_{n})^{+}(\sigma^{n})^{-})_{k+\frac{1}{2}}
-((\vartheta^{e}_{n})^{+}(\sigma^{n})^{+})_{k-\frac{1}{2}}\big)
-\sum_{k=1}^{K}\big((\pi^{e}_{n})^{-}(\tau^{n})^{-})_{k+\frac{1}{2}}
-((\pi^{e}_{n})^{-}(\tau^{n})^{+})_{k-\frac{1}{2}}\big).
\end{align}
\end{subequations}

\end{lem}
\textbf{Proof.} From the Galerkin orthogonality \eqref{1bch2k2}, we get
\begin{equation}\label{1bch2k}
\begin{split}
\bigg( &\sum_{j=1}^{S}W(\alpha_{j})\Delta\tau_{j}\delta_{t}^{\alpha_{j}}(\pi^{n}-\pi^{e}_{n}),\vartheta_{1}\bigg)_{D^{k}}
+\bigg(\sum_{j=1}^{S}W(\alpha_{j})\Delta\tau_{j}\delta_{t}^{\alpha_{j}}(\sigma-\sigma^{e}),\chi\bigg)_{D^{k}}
-\big(\Delta_{(\beta-2)/2}(\epsilon^{n}-\epsilon^{e}_{n}),\rho\big)_{D^{k}}
\\
&-\big(\Delta_{(\beta-2)/2}(\varphi^{n}-\varphi^{e}_{n}),\varrho\big)_{D^{k}}
+\big(\tau^{n}-\tau^{e}_{n},\phi_{x}\big)_{D^{k}} +\big(\sigma^{n}-\sigma^{e}_{n},\varphi_{x}\big)_{D^{k}}+\big(\vartheta^{n}-\vartheta^{e}_{n},\psi_{x}\big)_{D^{k}}
+\big(\pi^{n}-\pi^{e}_{n},\zeta_{x}\big)_{D^{k}}\\
&+\varepsilon_{2} \big(\sigma^{n}-\sigma^{e}_{n},\vartheta_{1}\big)_{D^{k}}- \varepsilon_{2}\big(\pi^{n}-\pi^{e}_{n},\chi\big)_{D^{k}}
+\big(\epsilon^{n}-\epsilon^{e}_{n},\phi\big)_{D^{k}}+\big(\tau^{n}-\tau^{e}_{n},\varphi\big)_{D^{k}}
+\big(\varpi^{n}-\varpi^{e}_{n},\varrho\big)_{D^{k}}\\
&+\big(\phi^{n}-\phi^{e}_{n},\rho\big)_{D^{k}}+\big(\varphi^{n}-\varphi^{e}_{n},\psi\big)_{D^{k}} +\big(\vartheta^{n}-\vartheta^{e}_{n},\zeta\big)_{D^{k}}+\varepsilon_{1}\big(\phi^{n}
-\phi^{e}_{n},\vartheta_{1}\big)_{D^{k}}-\varepsilon_{1}\big(\varpi^{n}-\varpi^{e}_{n},\chi\big)_{D^{k}}
\\
&-\sum_{k=1}^{K}\big(((\tau^{n}-\tau^{e}_{n})^{+}(\phi)^{-})_{k+\frac{1}{2}}
-((\tau^{n}-\tau^{e}_{n})^{+}(\phi)^{+})_{k-\frac{1}{2}}\big)
-\sum_{k=1}^{K}\big(((\sigma^{n}-\sigma^{e}_{n})^{-}(\varphi)^{-})_{k+\frac{1}{2}}
-((\sigma^{n}-\sigma^{e}_{n})^{-}(\varphi)^{+})_{k-\frac{1}{2}}\big)
\\
&-\sum_{k=1}^{K}\big(((\vartheta^{n}-\vartheta^{e}_{n})^{+}(\psi)^{-})_{k+\frac{1}{2}}
-((\vartheta^{n}-\vartheta^{e}_{n})^{+}(\psi)^{+})_{k-\frac{1}{2}}\big)-\sum_{k=1}^{K}
\big((\pi^{n}-\pi^{e}_{n})^{-}(\zeta)^{-})_{k+\frac{1}{2}}\\
&-((\pi^{n}-\pi^{e}_{n})^{-}(\zeta)^{+})_{k-\frac{1}{2}}\big)=0.\\
\end{split}
\end{equation}
We take the test functions
\begin{equation}\label{91h}
\begin{split}
\vartheta_{1}=\pi^{n},\quad\rho=\phi^{n}-\epsilon^{n},\quad \phi=\pi^{n},\quad \varphi=-\vartheta^{n},\quad \chi=\sigma^{n},\quad\varrho=\varpi^{n}-\varphi^{n},\quad \psi=-\sigma^{n},\quad \zeta=\tau^{n},
\end{split}
\end{equation}
 we obtain
\begin{equation}\label{1bch2k}
\begin{split}
\bigg( &\sum_{j=1}^{S}W(\alpha_{j})\Delta\tau_{j}\delta_{t}^{\alpha_{j}}(\pi^{n}-\pi^{e}_{n}),\pi^{n}\bigg)_{D^{k}}
+\bigg(\sum_{j=1}^{S}W(\alpha_{j})\Delta\tau_{j}\delta_{t}^{\alpha_{j}}(\sigma^{n}-\sigma^{e}_{n}),
\sigma^{n}\bigg)-\big(\Delta_{(\beta-2)/2}(\epsilon^{n}-\epsilon^{e}_{n}),\phi^{n}-\epsilon^{n}\big)_{D^{k}}
\\
&-\big(\Delta_{(\beta-2)/2}(\varphi^{n}-\varphi^{e}_{n}),\varpi^{n}-\varphi^{n}\big)_{D^{k}}
+\big(\tau^{n}-\tau^{e}_{n},\pi_{x}^{n}\big)_{D^{k}} -\big(\sigma^{n}-\sigma^{e}_{n},\vartheta_{x}^{n}\big)_{D^{k}}
-\big(\vartheta^{n}-\vartheta^{e}_{n},\sigma_{x}^{n}\big)_{D^{k}}
+\big(\pi^{n}-\pi^{e}_{n},\tau_{x}^{n}\big)_{D^{k}}\\
&+\varepsilon_{2} \big(\sigma^{n}-\sigma^{e}_{n},\pi\big)_{D^{k}}-\varepsilon_{2} \big(\pi_{n}-\pi^{e}_{n},\sigma\big)_{D^{k}}
+\big(\epsilon^{n}-\epsilon^{e}_{n},\pi\big)_{D^{k}}-\big(\tau^{n}-\tau^{e}_{n},\vartheta\big)_{D^{k}}
+\big(\varpi^{n}-\varpi^{e}_{n},\varpi^{n}-\varphi^{n}\big)_{D^{k}}\\
&+\big(\phi^{n}-\phi^{e}_{n},\phi^{n}-\epsilon^{n}\big)_{D^{k}}-\big(\varphi^{n}-\varphi^{e}_{n},\sigma^{n}\big)_{D^{k}} +\big(\vartheta^{n}-\vartheta^{e}_{n},\tau^{n}\big)_{D^{k}}+\varepsilon_{1}\big(\phi^{n}
-\phi^{e}_{n},\pi^{n}\big)_{D^{k}}-\varepsilon_{1}\big(\varpi^{n}-\varpi^{e}_{n},\sigma^{n}\big)_{D^{k}}
\\
&-\big(((\tau^{n}-\tau^{e}_{n})^{+}(\pi^{n})^{-})_{k+\frac{1}{2}}
-((\tau^{n}-\tau^{e}_{n})^{+}(\pi^{n})^{+})_{k-\frac{1}{2}}\big)+\big(((\sigma^{n}-\sigma^{e}_{n})^{-}(\vartheta^{n})^{-})_{k+\frac{1}{2}}
-((\sigma^{n}-\sigma^{e}_{n})^{-}(\vartheta^{n})^{+})_{k-\frac{1}{2}}\big)
\\
&+\big(((\vartheta^{n}-\vartheta^{e}_{n})^{+}(\sigma^{n})^{-})_{k+\frac{1}{2}}
-((\vartheta^{n}-\vartheta^{e}_{n})^{+}(\sigma^{n})^{+})_{k-\frac{1}{2}}\big)-\sum_{k=1}^{K}
\big((\pi^{n}-\pi^{e}_{n})^{-}(\tau^{n})^{-})_{k+\frac{1}{2}}
-((\pi^{n}-\pi^{e}_{n})^{-}(\tau^{n})^{+})_{k-\frac{1}{2}}\big)=0.\\
\end{split}
\end{equation}
Summing over $k$, simplify  by integration by parts and \eqref{fl}. This completes
the proof. $\quad\Box$

\begin{thm}\label{fgrhaaq}
Let $u$ be the exact solution of the problem \eqref{1bchj}, and let $u_{h}$ be the numerical solution of the fully discrete LDG scheme \eqref{1bch2}.  Then for small enough $h$, we have the following error estimates:
\begin{equation}\label{tt7mkcvw}
\begin{split}
&\|u(x,t_{n})-u_{h}^{n}\|_{L^{2}(\Omega)}\leq C\big(h^{N+1}+(\Delta t)^{1+\frac{\theta}{2}}+\theta^{2}\big).\\
\end{split}
\end{equation}
\end{thm}
\textbf{Proof}. We estimate the term $Q_{i}$, $\,\, i=1,...,4$. So we employ Young's inequality,  Lemma \ref{lga2} and the approximation results \eqref{sth}, we obtain
 \begin{equation}\label{tt12hcv}
\begin{split}
Q_{1}\leq &c_{5}\|\epsilon^{n}\|^{2}_{L^{2}(\Omega)}
+c_{6} \|\varphi^{n}\|^{2}_{L^{2}(\Omega)}+c_{1}\|\pi^{n}\|^{2}_{L^{2}(\Omega)}+
c_{2}\|\sigma^{n}\|^{2}_{L^{2}(\Omega)}+c_{3}\|\phi^{n}\|^{2}_{L^{2}(\Omega)}
+c_{4}\|\varpi^{n}\|^{2}_{L^{2}(\Omega)}.\\
\end{split}
\end{equation}
Using the definition of the numerical traces, \eqref{fl}, and the definitions of the projections $\mathcal{P}^{+},\mathcal{P}^{-}$ \eqref{prh}, we get
\begin{equation}\label{91chssf}
\begin{split}
Q_{2}=Q_{4}=0.
\end{split}
\end{equation}
From the approximation results \eqref{sth}, \eqref{sthzz} and \eqref{sthzzc} and Young's inequality, we obtain
\begin{equation}\label{tt12h}
\begin{split}
Q_{3}\leq &c_{5}\|\epsilon^{n}\|^{2}_{L^{2}(\Omega)}
+c_{6} \|\varphi^{n}\|^{2}_{L^{2}(\Omega)}+c_{1}\|\pi^{n}\|^{2}_{L^{2}(\Omega)}+
c_{2}\|\sigma^{n}\|^{2}_{L^{2}(\Omega)}\\
&+c_{3}\|\phi^{n}\|^{2}_{L^{2}(\Omega)}
+c_{4}\|\varpi^{n}\|^{2}_{L^{2}(\Omega)}+C\big(h^{2N+2}+(\Delta t)^{2+\theta}+\theta^{4}\big).\\
\end{split}
\end{equation}

Combining \eqref{tt12hcv}, \eqref{91chssf}, \eqref{tt12h} and \eqref{91chvv}, we obtain

\begin{equation}\label{91ch}
\begin{split}
&\bigg(\sum_{j=1}^{S}W(\alpha_{j})\Delta\tau_{j}\delta_{t}^{\alpha_{j}}\pi^{n},\pi^{n}\bigg)
+\bigg(\sum_{j=1}^{S}W(\alpha_{j})
\Delta\tau_{j}\delta_{t}^{\alpha_{j}}\sigma^{n},\sigma^{n}\bigg)
+\big(\Delta_{(\beta-2)/2}\epsilon^{n},\epsilon^{n}\big)
+\big(\Delta_{(\beta-2)/2}\varphi^{n},\varphi^{n}\big)\\
&\quad+\big(\phi^{n},\phi^{n}\big)
+\big(\varpi^{n},\varpi^{n}\big)\leq c_{5}\|\epsilon^{n}\|^{2}_{L^{2}(\Omega)}
+c_{6} \|\varphi^{n}\|^{2}_{L^{2}(\Omega)}+c_{1}\|\pi^{n}\|^{2}_{L^{2}(\Omega)}+
c_{2}\|\sigma^{n}\|^{2}_{L^{2}(\Omega)}+c_{3}\|\phi^{n}\|^{2}_{L^{2}(\Omega)}\\
&\quad\quad\qquad\qquad\qquad\qquad\qquad\qquad\qquad
+c_{4}\|\varpi^{n}\|^{2}_{L^{2}(\Omega)}+C\big(h^{2N+2}+(\Delta t)^{2+\theta}+\theta^{4}\big). \\
\end{split}
\end{equation}
Recalling Lemmas \ref{lg} and provided $c_{3}, c_{4}$ are sufficiently small such that $c_{3},c_{4}\leq1$, we obtain
\begin{equation}\label{91ch}
\begin{split}
\bigg(\sum_{j=1}^{S}W(\alpha_{j})\Delta\tau_{j}\delta_{t}^{\alpha_{j}}\pi^{n},\pi^{n}\bigg)
+\bigg(\sum_{j=1}^{S}W(\alpha_{j})
\Delta\tau_{j}\delta_{t}^{\alpha_{j}}\sigma^{n},\sigma^{n}\bigg)
\leq & c_{1}\|\pi^{n}\|^{2}_{L^{2}(\Omega)}+
 c_{2}\|\sigma^{n}\|^{2}_{L^{2}(\Omega)}\\
&+C\big(h^{2N+2}+(\Delta t)^{2+\theta}+\theta^{4}\big). \\
\end{split}
\end{equation}
It then follows that
\begin{equation}\label{1dv1}
\begin{split}
&\bigg(\sum_{j=1}^{S}
\frac{W(\alpha_{j})\Delta\tau_{j}}{\lambda_{j}}\pi^{n},\pi^{n}_{h}\bigg)+\bigg(\sum_{j=1}^{S}
\frac{W(\alpha_{j})\Delta\tau_{j}}{\lambda_{j}}\sigma^{n},\sigma^{n}_{h}\bigg)\\
&\leq
\bigg(\sum_{j=1}^{S}
\frac{W(\alpha_{j})\Delta\tau_{j}}{\lambda_{j}}\sum_{l=1}^{n-1}
(a_{n-l-1}^{\alpha_{j}}-a_{n-l}^{\alpha_{j}})\pi^{l},\pi^{n}\bigg)+\bigg(\sum_{j=1}^{S}
\frac{W(\alpha_{j})\Delta\tau_{j}}{\lambda_{j}}\sum_{l=1}^{n-1}
(a_{n-l-1}^{\alpha_{j}}-a_{n-l}^{\alpha_{j}})\sigma^{l},\sigma^{n}\bigg)\\
&\quad\quad+
\bigg(\sum_{j=1}^{S}
\frac{W(\alpha_{j})\Delta\tau_{j}}{\lambda_{j}}a_{n-1}^{\alpha_{j}}\sigma^{0},\sigma^{n}\bigg)+
\bigg(\sum_{j=1}^{S}
\frac{W(\alpha_{j})\Delta\tau_{j}}{\lambda_{j}}a_{n-1}^{\alpha_{j}}\pi^{0},\pi^{n}\bigg)+ c_{1}\|\pi^{n}\|^{2}_{L^{2}(\Omega)}+c_{2}\|\sigma^{n}\|^{2}_{L^{2}(\Omega)}\\
&\quad\quad+C\big(h^{2N+2}+(\Delta t)^{2+\theta}+\theta^{4}\big).
\end{split}
\end{equation}
Employing Young's inequality, we obtain
\begin{equation}\label{1dv1}
\begin{split}
&\|\pi^{n}\|^{2}_{L^{2}(\Omega)}+\|\sigma^{n}\|^{2}_{L^{2}(\Omega)}\\
&\leq \sum_{j=1}^{S}
\frac{W(\alpha_{j})\Delta\tau_{j}}{\lambda_{j}}Q\bigg(\sum_{l=1}^{n-1}
(a_{n-l-1}^{\alpha_{j}}-a_{n-l}^{\alpha_{j}})\big(\|\pi^{l}\|^{2}_{L^{2}(\Omega)}
+\|\sigma^{l}\|^{2}_{L^{2}(\Omega)}\big)\bigg)
+(cQ+\frac{1}{4})\sum_{j=1}^{S}
\frac{W(\alpha_{j})\Delta\tau_{j}}{\lambda_{j}}\big(\|\pi^{n}\|^{2}_{L^{2}(\Omega)}
\\
&\quad\quad+\|\sigma^{n}\|^{2}_{L^{2}(\Omega)}\big)+c\sum_{j=1}^{S}
\frac{W(\alpha_{j})\Delta\tau_{j}}{\lambda_{j}}Qa_{n-1}^{\alpha_{j}}h^{2N+2}+C\sum_{j=1}^{S}
\frac{W(\alpha_{j})\Delta\tau_{j}}{\lambda_{j}}Qa_{n-1}^{\alpha_{j}}\big(h^{2N+2}+(\Delta t)^{2+\theta}+\theta^{4}\big),\\
\end{split}
\end{equation}
  provided $c$ is sufficiently small such that $\frac{3}{4}-cQ>0$, we obtain that
\begin{equation}\label{1dv1cvmhj}
  \begin{split}
\|\pi^{n}\|^{2}_{L^{2}(\Omega)}+\|\sigma^{n}\|^{2}_{L^{2}(\Omega)}\leq &C\bigg(\sum_{j=1}^{S}
\frac{W(\alpha_{j})\Delta\tau_{j}}{\lambda_{j}}Q\bigg(\sum_{l=1}^{n-1}
(a_{n-l-1}^{\alpha_{j}}-a_{n-l}^{\alpha_{j}})\big(\|\pi^{l}\|^{2}_{L^{2}(\Omega)}
+\|\sigma^{l}\|^{2}_{L^{2}(\Omega)}\big)\bigg)\\
&+\sum_{j=1}^{S}
\frac{W(\alpha_{j})\Delta\tau_{j}}{\lambda_{j}}Qa_{n-1}^{\alpha_{j}}\big(h^{2N+2}+(\Delta t)^{2+\theta}+\theta^{4}\big)\bigg).\\
\end{split}
\end{equation}
Obviously the theorem holds for $n = 0$. Assume that it is valid for $n = 1, 2, ..., m-1$. Then,
by \eqref{1dv1cvmhj}, we have
\begin{equation}\label{1dv1}
  \begin{split}
\|\pi^{m}\|^{2}_{L^{2}(\Omega)}+\|\sigma^{m}\|^{2}_{L^{2}(\Omega)}&\leq C\bigg(\sum_{j=1}^{S}
\frac{W(\alpha_{j})\Delta\tau_{j}}{\lambda_{j}}Q\bigg(\sum_{l=1}^{m-1}
(a_{n-l-1}^{\alpha_{j}}-a_{n-l}^{\alpha_{j}})\big(\|\pi^{l}\|^{2}_{L^{2}(\Omega)}
+\|\sigma^{l}\|^{2}_{L^{2}(\Omega)}\big)\bigg)\\
&\qquad+\sum_{j=1}^{S}
\frac{W(\alpha_{j})\Delta\tau_{j}}{\lambda_{j}}Qa_{n-1}^{\alpha_{j}}\big(h^{2N+2}+(\Delta t)^{2+\theta}+\theta^{4}\big)\bigg)\\
&\leq C\bigg(\sum_{j=1}^{S}
\frac{W(\alpha_{j})\Delta\tau_{j}}{\lambda_{j}}Q\sum_{l=1}^{m-1}
(a_{n-l-1}^{\alpha_{j}}-a_{n-l}^{\alpha_{j}})(h^{2N+2}+(\Delta t)^{2+\sigma}+\sigma^{4})\\
&\qquad+\sum_{j=1}^{S}
\frac{W(\alpha_{j})\Delta\tau_{j}}{\lambda_{j}}Qa_{n-1}^{\alpha_{j}}\big(h^{2N+1}+(\Delta t)^{4+\theta}+\theta^{4}\big)\bigg)
\\
&=C\big(h^{2N+2}+(\Delta t)^{2+\theta}+\theta^{4}\big).
\end{split}
\end{equation}
Finally, by using triangle inequality and standard approximation theory, we can get \eqref{tt7mkcvw}. $\quad\Box$\\

\subsection{ LDG method for   the   coupled nonlinear
distributed-order time and space-fractional Schr\"{o}dinger equations}\label{s44}

In this section, we present and analyze the LDG method for the   coupled  nonlinear distributed-order time and space- fractional Schr\"{o}dinger equations

\begin{equation}\label{1bchcc}
\begin{split}
&i\mathcal{D}_{t}^{W(\alpha)}{}u_{1}- \varepsilon_{1}(-\Delta)^{\frac{\beta}{2}}u_{1}+ \varepsilon_{2}f(|u_{1}|^{2},|u_{2}|^{2})u_{1}=0,\\
&i\mathcal{D}_{t}^{W(\alpha)}{}u_{2}- \varepsilon_{3}(-\Delta)^{\frac{\beta}{2}}u_{2}+\varepsilon_{4} g(|u_{1}|^{2},|u_{2}|^{2})u_{2}=0.\\
\end{split}
\end{equation}
To define the local discontinuous Galerkin method, we rewrite  \eqref{1bchcc} as a first-order system:

\begin{equation}\label{1bch}
\begin{split}
&i\mathcal{D}_{t}^{W(\alpha)}{}u_{1}+\varepsilon_{1}e+ \varepsilon_{2}f(|u_{1}|^{2},|u_{2}|^{2})u_{1}=0,\\
&e=\Delta_{(\beta-2)/2}r,\quad r=\frac{\partial}{\partial x}s,\quad s=\frac{\partial}{\partial x}u_{1},\\
&i\mathcal{D}_{t}^{W(\alpha)}{}u_{2}+\varepsilon_{3}l+ \varepsilon_{4}g(|u_{1}|^{2},|u_{2}|^{2})u_{2}=0,\\
&l=\Delta_{(\beta-2)/2}w,\quad w=\frac{\partial}{\partial x}z,\quad z=\frac{\partial}{\partial x}u_{2}.\\
\end{split}
\end{equation}
We decompose the complex functions $u(x, t)$ and $v(x, t)$ into their real and imaginary parts. Setting
$u_{1}(x, t)=p(x, t) + iq(x, t)$ and $u_{2}(x, t) = \upsilon(x, t) + i\vartheta(x, t)$ in system \eqref{1bchcc}, we can obtain the following coupled system
\begin{equation}\label{1bx}
\begin{split}
&\mathcal{D}_{t}^{W(\alpha)}{}p+\varepsilon_{1}Q+ \varepsilon_{2}f(|u_{1}|^{2},|u_{2}|^{2})q=0,\\
&Q=\Delta_{(\beta-2)/2}r,\quad r=\frac{\partial}{\partial x}s,\quad s=\frac{\partial}{\partial x}q,\\
&\mathcal{D}_{t}^{W(\alpha)}{}q-\varepsilon_{1}H- \varepsilon_{2}f(|u_{1}|^{2},|u_{2}|^{2})p=0,\\
&H=\Delta_{(\beta-2)/2}w,\quad w=\frac{\partial}{\partial x}z,\quad z=\frac{\partial}{\partial x}p,\\
&\mathcal{D}_{t}^{W(\alpha)}{}\upsilon+\varepsilon_{3}L+ \varepsilon_{4}g(|u_{1}|^{2},|u_{2}|^{2})\vartheta=0,\\
&L=\Delta_{(\beta-2)/2}\rho,\quad\rho=\frac{\partial}{\partial x}\varpi,\quad \varpi=\frac{\partial}{\partial x}\vartheta,\\
&\mathcal{D}_{t}^{W(\alpha)}{}\vartheta-\varepsilon_{3}E- \varepsilon_{4}g(|u_{1}|^{2},|u_{2}|^{2})\upsilon=0,\\
&E=\Delta_{(\beta-2)/2}\xi,\quad\xi=\frac{\partial}{\partial x}\varrho,\quad \varrho=\frac{\partial}{\partial x}\upsilon.\\
\end{split}
\end{equation}

We define a fully discrete local discontinuous Galerkin scheme with as follows: find $p_{h}^{n}, q_{h}^{n}, Q^{n},r_{h}^{n},s_{h}^{n},H^{n}_{h},w_{h}^{n},z_{h}^{n}$,\\
$\upsilon_{h}^{n},\vartheta_{h}^{n},L^{n}_{h},\rho_{h}^{n},\varpi_{h}^{n},E_{h}^{n},\xi_{h}^{n}$,$\varrho_{h}^{n} \in V_{k}^{N}$, such that for all test functions $\vartheta_{1},\beta_{1},\phi,\varphi,\chi,\beta_{2},\psi$,\\
$\zeta,\gamma,\beta_{3},\delta,\varsigma,o,\beta_{4},\omega,\kappa\in V_{k}^{N}$,

\begin{equation}\label{1bch2ccf}
\begin{split}
&\bigg(\sum_{j=1}^{S}W(\alpha_{j})\Delta\tau_{j}\delta_{t}^{\alpha_{j}}p_{h}^{n},\vartheta_{1}\bigg)_{D^{k}}
+\varepsilon_{1}\big(Q_{h}^{n},\vartheta_{1}\big)_{D^{k}}+ \varepsilon_{2}\big(f(|u_{1h}^{n}|^{2},|u_{2h}^{n}|^{2})q_{h}^{n},\vartheta_{1}\big)_{D^{k}}=0,\\
&\big(Q_{h}^{n},\beta_{1}\big)_{D^{k}}=\big(\Delta_{(\beta-2)/2}r_{h}^{n},\beta_{1}\big)_{D^{k}},\\
&\big(r_{h}^{n},\phi\big)_{D^{k}}=-\big(s_{h}^{n},\phi_{x}\big)_{D^{k}}
+\big((\widehat{s}^{n}_{h}\phi^{-})_{k+\frac{1}{2}}
-(\widehat{s}^{n}_{h}\phi^{+})_{k-\frac{1}{2}}\big),\\
&\big(s_{h}^{n},\varphi\big)_{D^{k}}=-\big(q_{h}^{n},\varphi_{x}\big)_{D^{k}}
+\big((\widehat{q}^{n}_{h}\varphi^{-})_{k+\frac{1}{2}}
-(\widehat{r}^{n}_{h}\varphi^{+})_{k-\frac{1}{2}}\big),\\
&\bigg(\sum_{j=1}^{S}W(\alpha_{j})\Delta\tau_{j}\delta_{t}^{\alpha_{j}}q_{h}^{n},\chi\bigg)_{D^{k}}
-\varepsilon_{1}\big(H_{h}^{n},\chi\big)_{D^{k}}
- \varepsilon_{2}\big(f(|u_{1h}^{n}|^{2},|u_{2h}^{n}|^{2})p_{h}^{n},\chi\big)_{D^{k}}=0,\\
&\big(H_{h}^{n},\beta_{2}\big)_{D^{k}}=\big(\Delta_{(\beta-2)/2}w_{h}^{n},\beta_{2}\big)_{D^{k}},\\
&\big(w_{h}^{n},\psi\big)_{D^{k}}=-\big(z_{h}^{n},\psi_{x}\big)_{D^{k}}
+\big((\widehat{z}^{n}_{h}\psi^{-})_{k+\frac{1}{2}}
-(\widehat{z}^{n}_{h}\psi^{+})_{k-\frac{1}{2}}\big),\\
&\big(z_{h}^{n},\zeta\big)_{D^{k}}=-\big(p_{h}^{n},\zeta_{x}\big)_{D^{k}}
+\big((\widehat{p}^{n}_{h}\zeta^{-})_{k+\frac{1}{2}}
-(\widehat{p}^{n}_{h}\zeta^{+})_{k-\frac{1}{2}}\big),\\
&\bigg(\sum_{j=1}^{S}W(\alpha_{j})\Delta\tau_{j}\delta_{t}^{\alpha_{j}}\upsilon_{h}^{n},\gamma\bigg)_{D^{k}}
+\varepsilon_{3}\big(L_{h}^{n},\gamma\big)_{D^{k}}
+ \varepsilon_{4}\big(g(|u_{1h}^{n}|^{2},|u_{2h}^{n}|^{2})\vartheta_{h}^{n},\gamma\big)_{D^{k}}=0,\\
&\big(L_{h}^{n},\beta_{3}\big)_{D^{k}}=\big(\Delta_{(\beta-2)/2}\rho_{h}^{n},\beta_{3}\big)_{D^{k}},\\
&\big(\rho_{h}^{n},\delta\big)_{D^{k}}=-\big(\varpi_{h}^{n},\delta_{x}\big)_{D^{k}}
+\big((\widehat{\varpi}^{n}_{h}\delta^{-})_{k+\frac{1}{2}}
-(\widehat{\varpi}^{n}_{h}\delta^{+})_{k-\frac{1}{2}}\big),\\
&\big(\varpi_{h}^{n},\varsigma\big)_{D^{k}}=-\big(\vartheta_{h}^{n},\varsigma_{x}\big)_{D^{k}}
+\big((\widehat{\vartheta}^{n}_{h}\varsigma^{-})_{k+\frac{1}{2}}
-(\widehat{\vartheta}^{n}_{h}\varsigma^{+})_{k-\frac{1}{2}}\big),\\
&\bigg(\sum_{j=1}^{S}W(\alpha_{j})\Delta\tau_{j}\delta_{t}^{\alpha_{j}} \vartheta_{h}^{n},o\bigg)_{D^{k}}-\varepsilon_{3}\big(E_{h}^{n},o\big)_{D^{k}}
- \varepsilon_{4}\big(g(|u_{1h}^{n}|^{2},|u_{2h}^{n}|^{2})\upsilon_{h}^{n},o\big)_{D^{k}}=0,\\
&\big(E_{h}^{n},\beta_{4}\big)_{D^{k}}=\big(\Delta_{(\beta-2)/2}\xi_{h}^{n},\beta_{4}\big)_{D^{k}},\\
&\big(\xi_{h}^{n},\omega\big)_{D^{k}}=-\big(\varrho_{h}^{n},\omega_{x}\big)_{D^{k}}
+\big((\widehat{\varrho}^{n}_{h}\omega^{-})_{k+\frac{1}{2}}
-(\widehat{\varrho}^{n}_{h}\omega^{+})_{k-\frac{1}{2}}\big),\\
&\big(\varrho_{h}^{n},\kappa\big)_{D^{k}}=-\big(\upsilon_{h}^{n},\kappa_{x}\big)_{D^{k}}
+\big((\widehat{\upsilon}^{n}_{h}\kappa^{-})_{k+\frac{1}{2}}
-(\widehat{\upsilon}^{n}_{h}\kappa^{+})_{k-\frac{1}{2}}\big).\\
\end{split}
\end{equation}
The numerical traces $( \widehat{p}^{n}_{h}, \widehat{q}^{n}_{h}, \widehat{s}^{n}_{h},\widehat{z}^{n}_{h},\widehat{\upsilon}^{n}_{h},\widehat{\vartheta}^{n}_{h},\widehat{\varpi}^{n}_{h},\widehat{\varrho}^{n}_{h})$ are defined on interelement faces as the
alternating fluxes
\begin{equation}\label{flcf}
\begin{split}
&\widehat{p}^{n}_{h}=(p^{n}_{h})^{-},\,\,\widehat{s}^{n}_{h}
=(s^{n}_{h})^{+},\,\,\,\,
\widehat{q}^{n}_{h}=(q^{n}_{h})^{-},\,\,
\widehat{z}^{n}_{h}=(z^{n}_{h})^{+},\\
& \widehat{\upsilon}^{n}_{h}=(\upsilon^{n}_{h})^{-},\,\,\widehat{\varpi}^{n}_{h}
=(\varpi^{n}_{h})^{+},\,\,\,\,\widehat{\varrho}^{n}_{h}=(\varrho^{n}_{h})^{+},\,\,
\widehat{\vartheta}^{n}_{h}=(\vartheta^{n}_{h})^{-}.
\end{split}
\end{equation}
\begin{thm}\label{tt4ha} ($L^{2}$ stability). Suppose $u_{1}(x, t)=p(x, t) + iq(x, t)$ and $u_{2}(x, t) = \upsilon(x, t) + i\vartheta(x, t)$ and let  $u_{1h}^{n}, u_{2h}^{n}\in V_{k}^{N}$ be the approximation of $u_{1}(x,t_{n}),u_{2}(x,t_{n})$ then the solution to the scheme \eqref{1bch2ccf} and \eqref{flcf} satisfies the $L^{2}$ stability

$$
\|u_{1h}^{n}\|_{L^{2}(\Omega)}^{2}+\|u_{2h}^{n}\|_{L^{2}(\Omega)}^{2}
\leq C(\|u_{1h}^{0}\|_{L^{2}(\Omega)}^{2}+\|u_{2h}^{0}\|_{L^{2}(\Omega)}^{2}).$$

\end{thm}
\begin{thm}\label{fgrha}
Let $u_{1}(x,t_{n})$ and $u_{2}(x,t_{n})$  be the exact solutions of the linear coupled  fractional
Schr\"{o}dinger equations \eqref{1bchcc}, and let $u_{1h}^{n}$ and $u_{2h}^{n}$ be the numerical solutions of the fully discrete LDG scheme \eqref{1bch2ccf}.  Then for small enough $h$, we have the following error estimates:
\begin{equation}\label{tt7h}
\begin{split}
&\|u_{1}(x,t_{n})-u_{1h}^{n}\|_{L^{2}(\Omega)}+\|u_{2}(x,t_{n})-u_{2h}^{n}\|_{L^{2}(\Omega)}\leq C\big(h^{N+1}+(\Delta t)^{1+\frac{\theta}{2}}+\theta^{2}\big).\\
\end{split}
\end{equation}
\end{thm}
Theorem  \ref{fgrha} and \ref{tt4ha} can be proven by similar techniques as that in the proof of Theorem  \ref{tt4h} and \ref{fgrhaaq}. We will thus not
give the details here.
\section{Numerical examples}\label{s5}
In the following, we present  some numerical experiments to show the accuracy and the performance of the present LDG method for the distributed-order time and space-fractional convection-diffusion  and Schr\"{o}dinger type equations.
\begin{exmp}\label{ex1} Consider the distributed-order time and space-fractional diffusion equation
\begin{equation}\label{91gynb}
\begin{split}
&\mathcal{D}_{t}^{W(\alpha)}{}u(x,t)+\varepsilon(-\Delta)^{\frac{\beta}{2}}u(x,t)=g(x,t),\quad x\in[-1,1],\quad t\in(0,0.5],\\
&u(x,0) = 0,
\end{split}
\end{equation}
and the corresponding forcing term $g(x,t)$  is of the form

\begin{equation}\label{91}
\begin{split}
g(x,t)=\bigg((x^{2}-1)^{4}\mathcal{D}_{t}^{W(\alpha)}{}t^{2}+\varepsilon t^{2}(-\Delta)^{\frac{\beta}{2}} (x^{2}-1)^{4}\bigg),
\end{split}
\end{equation}
\end{exmp}
then the exact solution is $u(x,t)=t^{2}(x^{2}-1)^{4}$
with $\varepsilon=\frac{\Gamma(8-\beta)}{\Gamma(8)}$.\\
The problem is solved for several different values of $\beta$, polynomial orders $(N)$, and numbers of elements $(K)$. The errors and spatial convergence orders  are listed in Table \ref{Tab:ab1} and show that the LDG method can achieve the
accuracy of order $N +1$.

\begin{table}[!htb]
    \centering
\begin{center}
 \begin{tabular}{|c|| c c c c c c c |}
  \hline
 \hline
 &\multicolumn{7}{|c|}{$N=1$} \\

 \hline
  K&\multicolumn{7}{|c|}{
 \quad$5$\qquad\quad\qquad\qquad$10$\qquad\qquad\qquad\qquad$15$\qquad\qquad\qquad\qquad$20$ \qquad\quad\qquad\qquad}\\

  $\beta$ & $L^{2}$-Error  & $L^{2}$-Error &order & $L^{2}$-Error & order & $L^{2}$-Error &order\\ [0.5ex]

 \hline

1.2& 5.97e-02&8.6e-03 &2.8 &  3.4e-03& 2.29 & 1.8e-03& 2.21\\

1.4 & 2.84e-02&5.8e-03&2.28 & 2.5e-03& 2.08 & 1.3e-03 & 2.27
 \\
1.8 &1.91e-02 & 4.5e-03&2.09& 1.9e-03& 2.13 & 9.9e-04 & 2.27\\
%

 \hline
 \hline
 &\multicolumn{7}{|c|}{$N=2$} \\

 \hline
  K&\multicolumn{7}{|c|}{
 \quad$5$\qquad\quad \qquad\qquad $10$\qquad \qquad \qquad\qquad $15$\qquad\qquad\qquad\qquad $20$ \qquad\quad\qquad\qquad}\\

  $\beta$ & $L^{2}$-Error  & $L^{2}$-Error &order & $L^{2}$-Error & order & $L^{2}$-Error &order\\ [0.5ex]

 \hline

1.2 & 3.52e-02  &  4.3e-03& 3.03 & 1.2e-03& 3.15
&4.8e-04&3.19 \\

1.4 & 1.57e-02 & 2.1e-03& 2.9  & 5.9e-04 & 3.13
&2.6e-04&2.85 \\
1.8   & 1.45e-02& 1.8e-03& 3.01 &5.5e-04

 & 2.92&2.2e-04&3.19\\
 \hline
 \hline

\end{tabular}
\end{center}
\caption{ $L^{2}$-Error and order of convergence   for Example \ref{ex1} with $K$ elements and polynomial order $N$.}\label{Tab:ab1}
\end{table}

\begin{exmp}\label{ex1b}  We consider the distributed-order time and space-fractional Burgers' equation
\begin{equation}\label{91gyn2}
\begin{split}
&\mathcal{D}_{t}^{W(\alpha)}{}u(x,t)+\varepsilon(-\Delta)^{\frac{\beta}{2}}u(x,t)
+\frac{\partial}{\partial x}\bigg(\frac{u^{2}(x,t)}{2}\bigg)=g(x,t),\quad x\in[-1,1],\quad t\in(0,0.5],\\
&u(x,0) = 0,
\end{split}
\end{equation}
and the corresponding forcing term $g(x,t)$  is of the form
\begin{equation}\label{91}
\begin{split}
g(x,t)=\bigg((x^{2}-1)^4\mathcal{D}_{t}^{W(\alpha)}{}t^{2}+8t^{4}x(x^{2}-1)^7
+\varepsilon t^{2}(-\Delta)^{\frac{\beta}{2}}(x^{2}-1)^4\bigg).
\end{split}
\end{equation}
\end{exmp}
In this case, the exact solution will be  $u(x,t)=t^{2}(x^{2}-1)^4$ with  $\varepsilon=\frac{\Gamma(8-\beta)}{\Gamma(8)}$.\\
To complete the scheme, we choose a Lax-Friedrichs flux for the nonlinear term.  We take  $\Delta t=T/500,\,\theta=1/50$. The errors and  spatial convergence orders are listed in Table \ref{Tab:b}.
\begin{table}[!htb]
    \centering
\begin{center}
 \begin{tabular}{|c|| c c c c c c c |}
  \hline
 \hline
 &\multicolumn{7}{|c|}{$N=1$} \\

 \hline
  K&\multicolumn{7}{|c|}{
 \quad$10$\qquad\quad\qquad\qquad$20$\qquad\qquad\qquad\qquad$30$\qquad\qquad\qquad\qquad$40$ \qquad\quad\qquad\qquad}\\

  $\beta$ & $L^{2}$-Error  & $L^{2}$-Error &order & $L^{2}$-Error & order & $L^{2}$-Error &order\\ [0.5ex]

 \hline

1.2& 7.8e-03& 1.9e-03&2.04 &  8.5e-04& 1.98& 4.6e-04& 2.13\\

1.4 & 4.9e-03&1.1e-03&2.16  & 4.6e-04& 2.15 & 2.5e-04 & 2.12
 \\
1.8 &1.9e-03 & 5.1e-04&1.9& 2.2e-04& 2.07 & 1.2e-04 & 2.11\\
%

 \hline
 \hline
 &\multicolumn{7}{|c|}{$N=2$} \\

 \hline
  K&\multicolumn{7}{|c|}{
 \quad$10$\qquad\quad \qquad\qquad $20$\qquad \qquad \qquad\qquad $30$\qquad\qquad\qquad\qquad $40$ \qquad\quad\qquad\qquad}\\

  $\beta$ & $L^{2}$-Error  & $L^{2}$-Error &order & $L^{2}$-Error & order & $L^{2}$-Error &order\\ [0.5ex]

 \hline

1.2 & 3.4e-03  &  4.2e-04& 3.02 & 1.3e-04 & 2.89
&5.2e-05&3.19 \\

1.4 & 1.3e-03& 1.8e-04& 2.85  & 5.6e-05 & 2.88
&2.5e-05&2.8 \\
1.8   & 8.2e-04& 1.1e-04& 2.9 &3.1e-05

 &3.12&1.3e-08&3.02\\
 \hline
 \hline
\end{tabular}
\end{center}
\caption{ $L^{2}$-Error and order of convergence   for Example \ref{ex1b} with $K$ elements and polynomial order $N$.}\label{Tab:b}
\end{table}
Table  \ref{Tab:c10} provides some numerical results of the errors and the temporal convergence orders with $\beta = 1.2, 1.6$ respectively at $T = 0.5$ with $N=1$, $K=30$.  Numerical results of the errors and the numerical integration convergence orders in Table \ref{Tab:c20} with $\beta = 1.2, 1.6$ respectively at $T = 0.5$. From these tables, we can see that the convergence order of the scheme is $\mathcal{O}(h^{N+1}+(\Delta t)^{1+\frac{\theta}{2}}+\theta^{2})$, which matches the theoretical
convergence order when $\theta$ is small enough.

\begin{table}[!htb]
    \centering
\begin{center}
 \begin{tabular}{|c|| c c||c|| c c  |}
  \hline
 \hline

 \hline
  $\beta$&\multicolumn{5}{|c|}{   $\beta$=1.2\quad \quad  \quad  \qquad \quad $\beta$=1.6}\\
\hline\hline
  $\Delta t$ & $L^{2}$-Error & order & $\Delta t$&$L^{2}$-Error& order \\ [0.5ex]

 \hline
T/100 & 4.38e-03& -&       T/100&1.6e-03&- \\

T/200 & 2.21e-03 &0.99&   T/200&7.82e-04&1.03 \\
 T/400&1.1e-03&1.01&     T/400&4.1e-04&0.93\\

 \hline
 \hline

\end{tabular}
\end{center}
\caption{ $L^{2}$-Error and  temporal convergence orders for $u$ with $\beta = 1.2, 1.8$  at $T = 0.5$.}\label{Tab:c10}
\end{table}
\begin{table}[!htb]
    \centering
\begin{center}
 \begin{tabular}{|c|| c c||c|| c c  |}
  \hline
 \hline

 \hline
  $\beta$&\multicolumn{5}{|c|}{   $\beta$=1.2\quad \quad  \quad  \qquad \quad $\beta$=1.6}\\
\hline\hline
  $\theta$ & $L^{2}$-Error & order & $\theta$&$L^{2}$-Error& order \\ [0.5ex]

 \hline
1/10 &2.14e-02& -&       1/10&7.37e-03&- \\

1/20 &4.84e-03 &2.15&   1/20&2.01e-03&1.88 \\
 1/40&1.3e-03&1.9&     1/40&4.6e-04&2.13\\

 \hline
 \hline

\end{tabular}
\end{center}
\caption{ $L^{2}$-Error and   numerical integration convergence orders for $u$ with $\beta = 1.2, 1.8$  at $T = 0.5$.}\label{Tab:c20}
\end{table}
\begin{exmp}\label{ex1bbsch} Consider  the following nonlinear distributed-order time and space-fractional Schr\"{o}dinger equation
\begin{equation}\label{sch1}
\begin{split}
&i\mathcal{D}_{t}^{W(\alpha)}{}u(x,t)- \varepsilon(-\Delta)^{\frac{\beta}{2}}u+|u|^{2}u=g(x,t),\quad x\in[-1,1],\quad t\in(0,0.5],\\
&u(x,0)=0,
\end{split}
\end{equation}
and the corresponding forcing term $g(x,t)$  is of the form
\begin{equation}\label{91}
\begin{split}
g(x,t)=(1+i)\bigg(i(x^{2}-1)^5\mathcal{D}_{t}^{W(\alpha)}{}t^{2}-\varepsilon t^{2}(-\Delta)^{\frac{\beta}{2}}(x^{2}-1)^5
+2t^{6}(x^{2}-1)^{15}\bigg).
\end{split}
\end{equation}
\end{exmp}
The exact solution $u(x,t)=(1+i)t^{2}(x^{2}-1)^5$ with  $\varepsilon=\frac{\Gamma(10-\beta)}{\Gamma(10)}$. The errors and spatial convergence orders are listed in Table \ref{Tab:bsch}. 
\begin{table}[!htb]
    \centering
\begin{center}
 \begin{tabular}{|c|| c c c c c c c |}
  \hline
 \hline
 &\multicolumn{7}{|c|}{$N=1$} \\

 \hline
  K&\multicolumn{7}{|c|}{
 \quad$10$\qquad\quad\qquad\qquad$20$\qquad\qquad\qquad\qquad$30$\qquad\qquad\qquad\qquad$40$ \qquad\quad\qquad\qquad}\\

  $\beta$ & $L^{2}$-Error  & $L^{2}$-Error &order & $L^{2}$-Error & order & $L^{2}$-Error &order\\ [0.5ex]

 \hline

1.2& 1.23e-02& 4.61e-03&1.42 &  1.97e-03& 2.1& 1.1e-03&2.03\\

1.4 & 1.01e-02&2.51e-03&2.01  & 1.11e-03& 2.01 & 6.31e-04 &1.96
 \\
1.8 &7.31e-03 & 1.91e-03&1.94& 8.35e-04& 2.04 & 4.71e-04 & 1.99\\
%

 \hline
 \hline
 &\multicolumn{7}{|c|}{$N=2$} \\

 \hline
  K&\multicolumn{7}{|c|}{
 \quad$10$\qquad\quad \qquad\qquad $20$\qquad \qquad \qquad\qquad $30$\qquad\qquad\qquad\qquad $40$ \qquad\quad\qquad\qquad}\\

  $\beta$ & $L^{2}$-Error  & $L^{2}$-Error &order & $L^{2}$-Error & order & $L^{2}$-Error &order\\ [0.5ex]

 \hline

1.2 & 8.35e-03  & 1.21e-03& 2.79 & 3.55e-04 & 3.02
&1.41e-04&3.21 \\

1.4 & 6.24e-03& 9.23e-04& 2.76 & 2.79e-04 & 2.95
&1.13e-04&3.14\\
1.8   & 2.62e-03& 3.54e-04& 2.89 &1.13e-04

 &2.82&4.66e-05&3.08\\
 \hline
 \hline

\end{tabular}
\end{center}
\caption{ $L^{2}$-Error and order of convergence   for Example \ref{ex1bbsch} with $K$ elements and polynomial order $N$.}\label{Tab:bsch}
\end{table}
Table  \ref{Tab:c1} provides some numerical results of the errors and the temporal convergence orders with $\beta = 1.2, 1.6$ respectively at $T = 0.5$.  Numerical results of the
errors and the numerical integration convergence orders in Table \ref{Tab:c2} with $\beta = 1.2, 1.6$ respectively at $T = 0.5$. From these tables, we can see that the convergence
order of the scheme is $\mathcal{O}(h^{N+1}+(\Delta t)^{1+\frac{\theta}{2}}+\theta^{2})$, which matches the theoretical
convergence order when $\theta$ is small enough.


\begin{table}[!htb]
    \centering
\begin{center}
 \begin{tabular}{|c|| c c||c|| c c  |}
  \hline
 \hline

 \hline
  $\beta$&\multicolumn{5}{|c|}{   $\beta$=1.2\quad \quad  \quad  \qquad \quad $\beta$=1.6}\\
\hline\hline
  $\Delta t$ & $L^{2}$-Error & order & $\Delta t$&$L^{2}$-Error& order \\ [0.5ex]

 \hline
T/100 & 6.25e-03& -&       T/100&5.64e-03&- \\

T/200 & 3.12e-03 &1.00&   T/200&2.81e-03&1.01 \\
 T/400&1.58e-03&0.98&     T/400&1.25e-03&1.17\\

 \hline
 \hline

\end{tabular}
\end{center}
\caption{ $L^{2}$-Error and  temporal convergence orders for $u$ with $\beta = 1.2, 1.8$  at $T = 0.5$.}\label{Tab:c1}
\end{table}
\begin{table}[!htb]
    \centering
\begin{center}
 \begin{tabular}{|c|| c c||c|| c c  |}
  \hline
 \hline

 \hline
  $\beta$&\multicolumn{5}{|c|}{   $\beta$=1.2\quad \quad  \quad  \qquad \quad $\beta$=1.6}\\
\hline\hline
  $\theta$ & $L^{2}$-Error & order & $\theta$&$L^{2}$-Error& order \\ [0.5ex]

 \hline
1/10 &5.28e-02& -&       1/10&2.25e-02&- \\

1/20 &1.25e-02 &2.08&   1/20&5.4e-03&2.06 \\
 1/40&2.98e-03&2.07&     1/40&1.55e-03&1.80\\

 \hline
 \hline

\end{tabular}
\end{center}
\caption{ $L^{2}$-Error and   numerical integration convergence orders for $u$ with $\beta = 1.2, 1.8$  at $T = 0.5$.}\label{Tab:c2}
\end{table}

\begin{exmp}\label{ex6} We consider the   coupled nonlinear distributed-order time and space-fractional Schr\"{o}dinger equations
\begin{equation}\label{sch1}
\begin{split}
&i\mathcal{D}_{t}^{W(\alpha)}{}u_{1}(x,t)- \varepsilon_{1}(-\Delta)^{\frac{\beta}{2}}u_{1}(x,t)+2(|u_{1}(x,t)|^{2}+|u_{2}(x,t)|^{2})
u_{1}(x,t)=g_{1}(x,t),\,\, x\in[-1,1],\,\,t\in(0,0.5],\\
&i\mathcal{D}_{t}^{W(\alpha)}{}u_{2}(x,t)- \varepsilon_{2}(-\Delta)^{\frac{\beta}{2}}u_{2}(x,t)+4(|u_{1}(x,t)|^{2}
+|u_{2}(x,t)|^{2})u_{2}(x,t)=g_{2}(x,t),\, x\in[-1,1],\, t\in(0,0.5],\\
\end{split}
\end{equation}
 and the corresponding forcing terms $g_{1}(x,t)$ and $g_{2}(x,t)$  are of the form
\begin{equation}\label{91}
\begin{split}
&g_{1}(x,t)=(1+i)\bigg(i(x^{2}-1)^6\mathcal{D}_{t}^{W(\alpha)}{}t^{2}-
\varepsilon_{1}t^{2}(-\Delta)^{\frac{\beta}{2}}(x^{2}-1)^6+8t^{6}(x^{2}-1)^{18}\bigg),\\
&g_{2}(x,t)=(1+i)\bigg(i(x^{2}-1)^6\mathcal{D}_{t}^{W(\alpha)}{}t^{2}-\varepsilon_{2}(-\Delta)^{\frac{\beta}{2}}(x^{2}-1)^6
+16t^{6}(x^{2}-1)^{18}\bigg),
\end{split}
\end{equation}
\end{exmp}
to obtain an exact solutions $u_{1}(x,t)=(1+i)t^{2}(x^{2}-1)^6$ and $u_{2}(x,t)=(1+i)t^{2}(x^{2}-1)^6$
with $\beta=1.3,\,\varepsilon_{1}=\frac{\Gamma(13-\beta)}{2\Gamma(13)}$, $\varepsilon_{2}=\frac{\Gamma(13-\beta)}{2\Gamma(13)}$. The errors
and spatial convergence orders are listed in Tables \ref{Tab:e1} and \ref{Tab:e2}, confirming optimal $\mathcal{O}(h^{N+1})$ order of convergence across.
\begin{table}[!htb]
    \centering
\begin{center}
 \begin{tabular}{|c|| c c||c|| c c ||c||c c |}
  \hline
 \hline

 \hline
  N&\multicolumn{8}{|c|}{   N=1\quad \quad \quad\quad\quad \quad\qquad \qquad N=2\quad\quad\quad \qquad\quad \quad  \qquad \quad N=3}\\
\hline\hline
  K & $L^{2}$-Error & order & K &$L^{2}$-Error & order&K&$L^{2}$-Error& order \\ [0.5ex]

 \hline
10 & 4.23e-02& -&       10&1.45e-02&-&         10&7.87e-03&-  \\

20 & 9.98e-03 &2.08&   20&1.87e-03&2.96 &    20&4.65e-04&4.08 \\
 40&2.54e-03&1.97&     40&2.12e-04&3.14&    40&2.59e-05&4.17\\
 80&6.48e-04&1.97&     80&2.68e-05&2.98&    80&1.63e-06&3.99\\
 \hline
 \hline
\end{tabular}
\end{center}
\caption{ $L^{2}$-Error and order of convergence for $u_{1}$ with $K$ elements and polynomial order $N$.}\label{Tab:e1}
\end{table}

\begin{table}[!htb]
    \centering
\begin{center}
 \begin{tabular}{|c|| c c||c|| c c ||c||c c |}
  \hline
 \hline

 \hline
  N&\multicolumn{8}{|c|}{   N=1\quad \quad \quad\quad\quad \quad\qquad \qquad N=2\quad\quad\quad \qquad\quad \quad  \qquad \quad N=3}\\
\hline\hline
  K & $L^{2}$-Error & order & K &$L^{2}$-Error & order&K&$L^{2}$-Error& order \\ [0.5ex]

 \hline
10 & 3.98e-02 & -&       10&1.26e-02&-&    10&6.89e-03&-  \\

20 &9.19e-03 &2.12&   20&1.54e-03&3.03 &    20&3.84e-04&4.17\\
 40&2.23e-03&2.04&     40&1.82e-04&3.08&    40&2.39e-05&4.01\\
80&5.54e-04&2.01&     80&1.94e-05&3.23&    80&1.47e-06&4.02\\
 \hline
 \hline

\end{tabular}
\end{center}
\caption{ $L^{2}$-Error and order of convergence for $u_{2}$ with $K$ elements and polynomial order $N$.}\label{Tab:e2}
\end{table}
\section{Conclusions}\label{s6}
In this work, we developed and analyzed a local discontinuous Galerkin  method for solving the  distributed-order time and space-fractional convection-diffusion and Schr\"{o}dinger type equations
, and have proven the stability and  error estimates of these methods. Numerical experiments confirm that the optimal order of convergence is recovered. Future work will include the analysis of LDG method for two-dimensional fractional problems.
\\
\\


\end{document}